\documentclass{amsart}

\usepackage{amssymb}
\usepackage{dsfont}

\newtheorem{Lemma}{Lemma}
\newtheorem{Definition}[Lemma]{Definition}
\newtheorem{Theorem}[Lemma]{Theorem}
\newtheorem{Corollary}[Lemma]{Corollary}
\newtheorem{Proposition}[Lemma]{Proposition}

\newtheorem{Question}[Lemma]{Question}
\newtheorem{Remark}[Lemma]{Remark}
\newtheorem{Fact}[Lemma]{Fact}
\newtheorem{Claim}[Lemma]{Claim}

\newcommand{\ad}{\hbox{\rm ad}}
\newcommand{\Ad}{\hbox{\rm Ad}}
\newcommand{\eq}{\hbox{\rm eq}}

\newcommand{\SL}{\hbox{\rm SL}}
\newcommand{\SO}{\hbox{\rm SO}}

\newcommand{\Image}{\hbox{\rm Im}}

\renewcommand{\>}{\rangle}
\renewcommand{\o}{^{\circ}}

\newcommand{\C}{\mathds{C}}

\newcommand{\N}{\mathds{N}}
\newcommand{\R}{\mathds{R}}

\begin{document}

\title[Cartan subgroups]{Cartan subgroups of groups definable in o-minimal structures}

\author{El\'ias Baro}
\address{Departamento de \'Algebra, Facultad de Matem\'aticas, Universidad Complutense de Madrid, 
28040 Madrid, Spain}

\author{Eric Jaligot}
\address{Institut Fourier, CNRS, Universit\'e Grenoble I, 100 rue des maths, BP
74, 38402 St Martin d'H\`eres cedex, France}

\author{Margarita Otero}
\address{Departamento de Matem\'aticas, Universidad Aut\'onoma de Madrid, 28049
Madrid,
Spain}

\thanks{The first author is partially supported by MTM2011-22435 and Grupos UCM 910444, the third author is partially supported by MTM2011-22435 and PR2011-0048}




\date{\today}

\begin{abstract}
We prove that groups definable in o-minimal structures have Cartan subgroups, and only 
finitely many conjugacy classes 
of such subgroups. We also delineate with precision how these subgroups cover the ambient group. 
\end{abstract}

\maketitle 


\section{Introduction}\label{introduction}

If $G$ is an arbitrary group, a subgroup $Q$ of $G$ is called a {\em Cartan subgroup} 
({\em in the sense of Chevalley}) if it satisfies the two following conditions: 
\begin{itemize}
\item[$(1)$]
$Q$ is nilpotent and maximal with this property among subgroups of $G$. 
\item[$(2)$] 
For any subgroup $X\leq Q$ which is normal in $Q$ and of finite index in $Q$, 
the normalizer $N_G(X)$ of $X$ in $G$ contains $X$ as a finite index subgroup.   
\end{itemize}
The purely group-theoretic definition of a Cartan subgroup as above was designed by Chevalley 
in order to capture critical properties of very specific subgroups of Lie groups. 

In connected reductive algebraic groups over algebraically closed fields and in connected compact 
real Lie groups, Cartan subgroups correspond typically to centralizers of maximal tori 
and it is well known that they are connected. 
It is however worth emphasizing at the outset that in real Lie groups 
Cartan subgroups need not be connected in general, a point also noticed by Chevalley 
in the introduction of \cite[Chapitre VI]{MR0068552}: 
``{\em Il convient de noter que les groupes de Cartan de $G$ ne sont en g\'en\'eral pas connexes.}" The diagonal subgroup of $\SL_2(\R)$ is maybe the first example of a nonconnected Cartan 
subgroup that one should bear in mind. Most of the difficulties for the study of 
these subgroups in the past, notably in the early work of Cartan, have been 
this failure of connectedness. This is something that will eventually need considerable 
attention in the present paper as well. 

We are going to study Cartan subgroups from the model-theoretic point of view of 
groups definable in an o-minimal structure, that is 
a first-order structure $\mathcal M=\<M,\leq,\cdots\>$ equipped with a total, dense, 
and without end-points definable order $\leq$ and such that every definable subset of $M$ is a 
boolean combination of intervals with end-points in $M\cup\{\pm\infty\}$. 
The most typical example of an o-minimal structure is of course the ordered field $\R$ of the reals, 
but there are richer o-minimal structures, such as the field of the reals equipped in addition 
with the exponential function \cite{MR1398816}. 

In order to deal with the non-connectedness of Cartan subgroups in general, 
we will use the following notion. 
If $G$ is a group definable in an arbitrary structure $\mathcal M$, then we say that it is 
{\em definably connected} 
if and only if it has no proper subgroup of finite index definable in the sense of $\mathcal M$. 
Now, a subgroup of a group $G$ definable in $\mathcal M$ is called a {\em Carter} subgroup of $G$ 
if it is definable and definably connected (in the sense of $\mathcal M$ as usual), 
and nilpotent and of finite index in its normalizer in $G$. All the notions 
of definability depend on a ground structure $\mathcal M$, 
which in the present paper will typically be an o-minimal structure. 
The notion of a Carter subgroup first appeared in the case of finite groups as 
{\em nilpotent and selfnormalizing} subgroups. A key feature is that, in the case of 
finite solvable groups, they exist and are conjugate \cite{Carter61/62}. 
For infinite groups, the notion we are adopting here, incorporating definability 
and definable connectedness, comes from the theory of groups of finite Morley rank. 
That theory is another classical branch of group theory in model theory, 
particularly designed at generalizing algebraic groups 
over algebraically closed fields. We note that the selfnormalization from the finite case 
becomes an almost selfnormalization property, and indeed the finite group $N_G(Q)/Q$ 
associated to a Carter subgroup $Q$ typically generalizes the notion of the {\em Weyl group} 
relative to $Q$. This is something that will also make perfect sense here in the case of groups definable in 
o-minimal structures. 

We will see shortly in Section \ref{SectionCarterCartan} 
that for groups definable in o-minimal structures, and actually for groups with the mere 
descending chain condition on definable subgroups, there is an optimal 
correspondence between Cartan subgroups and Carter subgroups: 
the latter ones are exactly the definably connected 
components of the former ones. 
In particular Cartan subgroups are automatically definable subgroups, a point not following from 
the definition of Chevalley in general, but which is always going to be true here. 

In Sections \ref{DimensionUnions}-\ref{SectionWeakGenNilp} we will 
relate Cartan and Carter subgroups to a well behaved notion of dimension 
for sets definable in an o-minimal structure, notably to 
{\em weak genericity} (having maximal dimension) 
or to {\em largeness} (having smaller codimension). 
We will mainly develop their {\em generous} analogs, where one actually considers 
the weak genericity or the largeness of the union of conjugates of a given set. 
The technics and results here will be substantial adaptations and generalizations from 
\cite{MR2225896, MR2537672} in the finite Morley rank case, and our arguments for 
Cartan and Carter subgroups of groups definable in o-minimal structure will highly depend on 
dimensional computations and generosity arguments. 
We will make such dimensional computations in a rather axiomatic 
framework, essentially with the mere existence of a definable and additive dimension, 
since they apply as such in many other contexts 
(groups of finite Morley rank, groups in supersimple theories of 
finite rank, groups definable over the $p$-adics...). 

Our main result can be summarized as follows. 

\begin{Theorem}\label{MainTheo}
Let $G$ be a group definable in an o-minimal structure. Then Cartan subgroups of $G$ exist, 
are definable, and fall into finitely many conjugacy classes. 
\end{Theorem}

Our proof of Theorem \ref{MainTheo} will also strongly depend on the main structural theorem about groups 
definable in o-minimal structures. It says in essence that any definably connected group 
$G$ definable in an o-minimal structure is, 
modulo a largest normal solvable (and definable) subgroup $R(G)$, a direct product of 
finitely many definably simple groups which are essentially ``known" as groups of Lie type. 
Hence our proof will consist in an analysis of the interplay between these definably simple factors 
and the relevant definably connected solvable subgroups of $G$. Results specific about groups 
definable in an o-minimal structure which are used here will be reviewed in 
Section \ref{SectionFactsGpsInoMinStruct}. 

A large part of the work will thus be concerned with the case of definably connected 
solvable groups. In this case we will make a strong use of the previously mentioned 
largeness and generosity arguments. Mixing them with more algebraic inductive 
arguments inspired by \cite{MR1765797} in the finite Morley rank case, we will obtain 
the following result in Section \ref{SectionSolvableCase}. 

\bigskip
\noindent
{\bf Theorem \ref{MainTheoSolvGps}}
{\em 
Let $G$ be a definably connected solvable group definable in an o-minimal structure. 
Then Cartan subgroups of $G$ exist and are conjugate, and they are definably connected 
and selfnormalizing. 
Moreover, they are largely generous in the following strong sense: for any Cartan 
subgroup $Q$, the (definable) set of elements of $Q$ contained in a unique conjugate of 
$Q$ is large in $Q$ and largely generous in $G$.  
}

\bigskip
A definably connected group is {\em semisimple} if it has a finite center and 
modulo that center abelian normal subgroups are trivial. Semisimplicity is a first-order 
property, and the main theorem about groups definable in o-minimal 
structures actually says that any such semisimple group with a trivial center 
is a direct product of definably simple groups, with each factor a ``known" group of Lie type modulo 
certain elementary equivalences. We will review certain facts more or less classical 
about Cartan subgroups of Lie groups in Section \ref{SectionLieGroups}. 
In Section \ref{SectionFromLieToDef} we will transfer the theory of Cartan subgroups of Lie groups 
to definably simple groups and get a quite complete description of Cartan subgroups 
of definably simple groups definable in o-minimal structures. 

In Section \ref{SectionSemiSimpleCase} we will elaborate further on the definably simple case 
to get a similarly quite complete description of Cartan subgroups of semisimple groups 
definable in o-minimal structures, obtaining the following general theorem. 

\bigskip
\noindent
{\bf Theorem \ref{transfer2} (lite)}
{\em 
Let $G$ be a definably connected semisimple group definable in an o-minimal structure. 
Then $G$ has definable Cartan subgroups and the following holds.
\begin{itemize}
\item[$(1)$] 
$G$ has only finitely many conjugacy classes of Cartan subgroups.
\item[$(2)$] 
If $Q_1$ and $Q_2$ are Cartan subgroups and  $Q_1\o=Q_2\o$, then  $Q_1=Q_2$.
\item[$(3)$] 
If $Q$ is a Cartan subgroup, then $Z(G)\leq Q$, $Q'\leq Z(G)$, and $Q\o\leq Z(Q)$. 
\item[$(4)$] 
If $Q$ is a Cartan subgroup and $a\in Q$, then $aQ\o$ is weakly generous. 
\item[$(5)$] 
The union of all Cartan subgroups, which is definable by {\em(1)}, is large. 
\end{itemize}
}

\bigskip
The general case of a definably connected group $G$ definable in an o-minimal structure 
will be considered in Section \ref{SectionGeneralCase}. In this case we have both $G$ not solvable 
and not semisimple, or in other words 
$$G/R\o(G)\neq 1 \mbox{~and~}R\o(G)\neq 1.$$
In that case Theorem \ref{MainTheo} follows rapidly from Theorems \ref{MainTheoSolvGps} 
and \ref{transfer2}, but some natural questions will remain without answer here. The most important one 
is maybe the following: if $Q$ is a Cartan subgroup of $G$, is it the case that 
$QR\o(G)/R\o(G)$ is a Cartan subgroup of the semisimple quotient $G/R\o(G)$? 
This question is indeed equivalent to the fact that Cartan subgroups of $G/R\o(G)$ 
are exactly of the form $QR\o(G)/R\o(G)$ for some Cartan subgroup $Q$ of $G$. 
We will only manage to prove that for a Cartan subgroup $Q$ of $G$, the group 
$QR\o(G)/R\o(G)$ is a finite index subgroup of a Cartan subgroup of $G/R\o(G)$, obtaining 
in particular the expected lifting for the corresponding Carter subgroups. 
Getting the exact lifting of Cartan subgroups seems to be related to interesting new problems of 
representation theory in a definable context. 
In any case, we will mention all what we managed to prove on the correlations between Cartan 
subgroup of $G$ and of $G/R\o(G)$, trying also to work with a not necessarily definably 
connected ambient group $G$ when possible. We will conclude in Section \ref{SectionFinalComments} with further 
comments on certain specialized topics, including algebraic or compact factors, 
Weyl groups relative to the various Cartan subgroups, and parameters. 

In this paper definability always means definability with parameters. 
We refer to \cite{MR2436142} for a complete introduction to groups definable in 
o-minimal structures. We insist that everything is done here for groups {\em definable} 
(as opposed to interpretable) in an arbitrary o-minimal structure. This is because the theory of 
groups in o-minimal structure has been developed in this slightly restricted context since 
\cite{MR961362}, where it is shown that definable groups can be equipped with a nice definable 
manifold structure making them topological groups. An arbitrary o-minimal structure does 
not eliminate imaginaries in general, but any group definable in an arbitrary 
o-minimal structure eliminates imaginaries, and actually has definable choice functions in a very 
strong sense \cite[Theorem 7.2]{MR2006422}. 
In particular, imaginaries coming from a group definable in an o-minimal structure will 
always be considered as definable in the sequel, and can be equipped with a finite dimension as any 
definable set. 
We refer to \cite[Chapter 4]{MR1633348} or \cite{MR961362} for the dimension of sets definable 
in o-minimal structures. 

Since we already gave the organization of the paper, let us immediately enter into its core. 

\section{Cartan subgroups and Carter subgroups}\label{SectionCarterCartan}

We first consider the relations between Cartan and Carter subgroups of groups 
definable in o-minimal structures. Actually, by \cite[Remark 2.13]{MR961362}, 
such groups satisfy the {\em descending chain condition} on definable subgroups 
($dcc$ for short), and we will analyze these relations in the more natural context of groups 
with the $dcc$. 
Throughout the present section, $G$ is a group 
definable in a structure $\mathcal M$ and definability may refer to $\mathcal{M}^{\eq}$, and we say that it 
satisfies the $dcc$ if any strictly descending chain of definable subgroups is stationary after 
finitely many steps. Notice that the $dcc$ always pass to quotients by definable normal subgroups. 

We first list some general facts needed in the sequel. 

\begin{Fact}\label{ActionConOnFinite}
{\bf \cite[Fact 3.1]{BaroJalOteroCommutators}}
Let $G$ be a definably connected group. 
\begin{itemize}
\item[$(a)$]
Any definable action of $G$ on a finite set is trivial. 
\item[$(b)$] 
If $Z(G)$ is finite, then $G/Z(G)$ is centerless. 
\end{itemize}
\end{Fact}

In a group with the $dcc$, any subset $X$ is contained in a smallest definable subgroup 
$H(X)$ called the {\em definable hull} of $X$: take $H(X)$ to be the intersection of 
all definable subgroups of $G$ containing $X$. 

\begin{Fact}\label{FactGenGdcc}
{\bf \cite[3.3 \& 3.4]{BaroJalOteroCommutators}}
Let $G$ be a group with the dcc and $X$ a subset of $G$. 
\begin{itemize}
\item[$(a)$]
If $X$ is $K$-invariant for some subset $K$ of $G$, then $H(X)$ is $K$-invariant as well.  
\item[$(b)$]
If $X$ is a nilpotent subgroup of $G$, then $H(X)$ is nilpotent of the same nilpotency class. 
\end{itemize}
\end{Fact}
We now mention an infinite version of the classical {\em normalizer condition} in finite nilpotent groups. 

\begin{Lemma}\label{NormalizerCondition}
Let $G$ be a nilpotent group with the $dcc$ on definable subgroups, or merely such that each definable subgroup 
has a definably connected definable subgroup of finite index. 
If $H$ is a definable subgroup of infinite index in $G$, then $N_G(H)/H$ is infinite. 
\end{Lemma}
\proof
For instance, one may argue formally as in \cite[Proposition 1.12]{MR1827833}. 
\qed

\begin{Lemma}\label{Carter-Cartan-Gen}
Let $G$ be a group with the $dcc$. 
\begin{itemize}
\item[$(a)$]
If $Q$ is a maximal nilpotent subgroup of $G$, then $Q$ is definable. 
\item[$(a')$]
If $Q$ is a Cartan subgroup of $G$, then $Q$ is definable and $Q\o$ is a Carter subgroup of $G$. 
\item[$(b)$]
If $Q$ is a Carter subgroup of $G$, then $Q$ is contained in a maximal nilpotent 
subgroup $\tilde{Q}$ of $G$, and any such subgroup $\tilde{Q}$ is a Cartan subgroup of $G$ with $[\tilde{Q}]\o=Q$. 
\end{itemize}
\end{Lemma}
\proof
$(a)$. By Fact \ref{FactGenGdcc}$(b)$. 

$(a')$. $Q$ is definable by item $(a)$. Since $Q\o$ is a normal subgroup of $Q$ 
of finite index in $Q$, $Q\o$ is a finite index subgroup of $N_G(Q\o)$, and 
$Q\o$ is a Carter subgroup of $G$. 

$(b)$. A definable nilpotent subgroup $H$ containing $Q$ must satisfy $H\o=Q$ by Lemma \ref{NormalizerCondition}, and 
thus $H\leq N_G(H\o)=N_G(Q)$. Now Fact \ref{FactGenGdcc}$(b)$ implies that any nilpotent subgroup $H$ containing $Q$ satisfies 
$Q\leq H \leq N_G(Q)$. Since $N_G(Q)/Q$ is finite, there are maximal such subgroups, proving our first claim. 

Now fix any such maximal nilpotent subgroup $\tilde Q$. It is definable by item $(a)$ and we have already seen that 
$Q={[\tilde Q]\o}$, and $\tilde{Q} \leq N_G([\tilde Q]\o)=N_G(Q)$. 
We now check that $\tilde Q$ is a Cartan subgroup. 
Let $X$ be any normal subgroup of finite index of $\tilde Q$. 
We first  observe that $H\o(X)=Q$: since $\tilde Q$ is definable we get 
$H\o(X)\leq{[\tilde Q]\o}=Q$, and since $H\o(X)$ must have finite index in $\tilde Q$ 
we get the desired equality. Now by Fact \ref{FactGenGdcc}$(a)$ 
$N_G(X)$ normalizes $H\o(X)=Q$, so $X\leq N_G(X)\leq N_G(Q)$. Since 
$X$ has finite index in $\tilde Q$ and $\tilde Q$ has finite 
index in $N_G(Q)$, $X$ has finite index in $N_G(Q)$, and in particular 
$X$ has finite index in $N_G(X)$. 
\qed

Applying Lemma \ref{Carter-Cartan-Gen}, we have thus that in 
groups definable in o-minimal structures Carter subgroups are {\em exactly} 
the definably connected components of Cartan subgroups, with the latter ones always definable. 
We also note that Lemma \ref{Carter-Cartan-Gen}$(a)$ gives the automatic definability of 
unipotent subgroups in many contexts, but that such unipotent subgroups 
are in general not almost selfnormalizing. We also note that if $Q$ is a maximal nilpotent subgroup, 
then it is a Cartan subgroup if and only if $Q\o$ is a Carter subgroup, by 
Lemma \ref{Carter-Cartan-Gen}. Finally, a selfnormalizing Carter subgroup must be a Cartan 
subgroup by Lemma \ref{Carter-Cartan-Gen}$(b)$, and a definably connected Cartan 
subgroup must be a Carter subgroup. 

Definably connected nilpotent groups definable in o-minimal structures are 
divisible by \cite[Theorem 6.10]{MR2006422}, so it is worth bearing in mind that the following always applies in 
groups definable in o-minimal structures. 

\begin{Fact}\label{StructNilpGps}
{\bf \cite[Lemma 3.10]{BaroJalOteroCommutators}}
Let $G$ be a nilpotent group with the $dcc$ and such that $G\o$ is divisible. 
Then $G=B\ast G\o$ (central product) for some finite subgroup $B$ of $G$. 
\end{Fact}

When Fact \ref{StructNilpGps} applies, one can strengthen Lemma \ref{Carter-Cartan-Gen}$(b)$ as follows. 
Again the following statement is valid in groups definable in an o-minimal structure, because they cannot contain 
an infinite increasing chain of definably connected subgroups (by the existence of a well behaved notion of dimension 
\cite[Corollary 2.4]{MR2436142}). 

\begin{Lemma}\label{existmaxnil}
Let $G$ be a group with the $dcc$. Assume that definably connected definable nilpotent 
subgroups of $G$ are divisible, and that $G$ contains no infinite increasing chain of 
such subgroups. 
Then any definably connected definable nilpotent subgroup of $G$ is contained in a maximal 
nilpotent subgroup of $G$. 
\end{Lemma}
\proof
Let $N$ be a definably connected nilpotent subgroup of $G$. By assumption, 
$N$ is contained in a definably connected definable nilpotent subgroup $N_1$ which 
is maximal for inclusion. It suffices to show that $N_1$ is then contained in a maximal nilpotent 
subgroup of $G$, and by Fact \ref{FactGenGdcc}$(b)$ 
we may consider only definable nilpotent subgroups 
containing $N_1$. It suffices then to show that any strictly increasing chain of definable nilpotent 
subgroups $N_1< N_2 < \cdots$ is stationary after finitely many steps. 

Assume towards a contradiction that $N_1< N_2< \cdots$ is such an infinite increasing 
chain of definable nilpotent subgroups. Recall that $N_1=N_1\o$, and notice also 
that $N_i\o=N_1$ for each $i$, since $N_1$ is maximal subject to being definably connected and containing $N$. 
By Fact \ref{StructNilpGps}, each $N_i$ has the form 
$B_i \ast N_1$ for some finite subgroup $B_i \leq N_i$, and in particular 
$N_i\leq {C_G(N_1)\cdot N_1}$. We may thus replace $G$ by the definable subgroup 
${C_G(N_1)\cdot N_1}$. 

Let $X$ be the union of the groups $N_i$. 
Working modulo the normal subgroup $N_1$, we have an increasing 
chain of finite nilpotent groups. 
Now $X/{N_1}$ is a periodic locally nilpotent group with the $dcc$ on centralizers, and by 
\cite[Theorem A]{MR549936} it is nilpotent-by-finite. 
Replacing $X$ by a finite index subgroup of $X$ if necessary, we may thus assume 
$X/{N_1}$ nilpotent and infinite. Since $G={C_G(N_1)\cdot N_1}$, 
the nilpotency of $X/{N_1}$ and of $N_1$ forces $X$ to be nilpotent (of nilpotency class bounded 
by the sum of that of $X/N_1$ and $N_1$). Replacing $X$ by $H(X)$, 
we may now assume with Fact \ref{FactGenGdcc}$(b)$ that $X$ is a 
definable nilpotent subgroup containing $N_1$ as a 
subgroup of infinite index. Then $N_1<X\o$, a contradiction to the maximality of $N_1$. 
\qed

Before moving ahead, it is worth mentioning concrete examples of Cartan subgroups of 
real Lie groups to be kept in mind in the present paper. 
In $\SL_2(\R)$ there are up to conjugacy two Cartan subgroups, the subgroup of diagonal 
matrices $Q_1\simeq \R^{\times}$, noncompact and not  connected, with 
corresponding Carter subgroup $Q_1\o\simeq \R^{>0}$, and $Q_2=\SO_2(\R)$ 
isomorphic to the circle group, 
compact and connected and hence also a Carter subgroup. 
More generally, and referring to \cite[p.141-142]{MR568880} for more details, 
the group $\SL_n(\R)$ has up to conjugacy $\left[\frac{n}{2}\right]+1$ Cartan subgroups 
$$Q_j\simeq[\C^{\times}]^{j-1}\!\times[\R^{\times}]^{n-2j+1}\mbox{~where~}1\leq j\leq {\left[\frac{n}{2}\right]+1},$$ 
unless $Q_{\frac{n}{2}+1}\simeq[\C^{\times}]^{\frac{n}{2}-1}\!\times {\SO_2(\R)}$ if $n=2(j-1)$. 

We will need the following lemma relating the center to Cartan and Carter subgroups. 
For any group $G$ we define the iterated centers $Z_n(G)$ as follows: $Z_0(G)=\{1\}$ and 
by induction $Z_{n+1}(G)$ is the preimage in $G$ of the center $Z(G/Z_n(G))$ of $G/Z_n(G)$. 

\begin{Lemma}\label{LemElemCartanModZn}
Let $G$ be a group and for $n\geq 0$ let $Z_n :=Z_n(G)$. 
\begin{itemize}
\item[$(a)$]
If $Q$ is a Cartan subgroup of $G$, then 
$Z_n\leq Q$ and $Q/Z_n$ is a Cartan subgroup of $G/Z_n$, 
and conversely every Cartan subgroup 
of $G/Z_n$ has this form. 
\item[$(b)$]
If $G$ satisfies the dcc, 
then Carter subgroups of $G/Z_n$ are exactly subgroups of the form $Q\o Z_n/Z_n$, 
for $Q$ a Cartan subgroup of $G$. 
\end{itemize}
\end{Lemma}
\proof
We may freely use the fact that the
preimage in $G$ of a nilpotent subgroup of $G/Z_n$ is nilpotent. 

$(a)$. Clearly $Z_n\leq Q$ by maximal nilpotence of $Q$. Clearly also, $Q/Z_n$ 
is nilpotent maximal in the quotient 
$\overline{G}=G/Z_n$. Let $\overline{X}$ be a normal subgroup of finite index of 
$\overline{Q}=Q/Z_n$, for some subgroup $X$ of $G$ containing $Z_n$. 
The preimage in $G$ of $N_{\overline{G}}(\overline{X})$ normalizes 
$X$, which clearly is normal and has finite index in $Q$. Since $Q$ is a Cartan subgroup of $G$, 
we easily get that $\overline{X}$ has finite index in $N_{\overline{G}}(\overline{X})$. 

Conversely, let $Q$ be a subgroup of $G$ containing $Z_n$ such that $Q/Z_n$ is a 
Cartan subgroup of $\overline{G}=G/Z_n$. Clearly $Q$ has to be 
maximal nilpotent in $G$. Let $X$ be a normal finite index subgroup of $Q$. $N_G(X)$ normalizes $\overline{X}$ 
modulo $Z_n$, so it must contain $\overline{X}$ as a finite index subgroup, 
and then $X$ is also a finite index subgroup of $N_G(X)$. 

$(b)$. By item $(a)$ Cartan subgroups of $G/Z_n$ are exactly of the form $Q/Z_n$ for a Cartan subgroup 
$Q$ of $G$ containing $Z_n$. So Carter subgroups of $G/Z_n$ are by Lemma \ref{Carter-Cartan-Gen} 
exactly of the form $[Q/Z_n]\o=Q\o Z_n/Z_n$, for $Q$ a Cartan subgroup of $G$. 
\qed

Finally, we will also use the following lemma describing Cartan subgroups of central products. 

\begin{Lemma}\label{LemElemCartanCentProd}
Let $G=G_1\ast \cdots \ast G_n$ be a central product of finitely many and pairwise commuting groups $G_i$. 
Then Cartan subgroups of $G$ are exactly of the form 
$Q_1\ast \cdots \ast Q_n$ where each $Q_i$ is a Cartan subgroup of $G_i$. 
\end{Lemma}
\proof
It suffices to prove our claim for $n=2$. For $i=1$ and $2$ and $X$ an arbitrary subset of $G$, let 
$\pi_i(X)=\{{g\in G_i}~|~\exists{h\in G_{i+1}}~gh\in X\}$,  
where the indices $i$ are of course considered modulo $2$.  
It is clear that when $X$ is a subgroup of $G$, $\pi_i(X)$ is a subgroup $G_i$. If $X$ is nilpotent (of nilpotency 
class $k$), then $\pi_i(X)$ is nilpotent (of nilpotency class at most $k+1$): it suffices to 
consider $G/G_{i+1}$ and to use the fact that ${G_1 \cap G_2}\leq Z(G_i)$. 

Let $Q$ be a Cartan subgroup of $G_1\ast G_2$. Since 
$Q\leq {\pi_1(Q) \ast \pi_2(Q)}$, the maximal nilpotence of $Q$ forces equality. 
Now it is clear that each $\pi_i(Q)$ is maximal nilpotent in $G_i$, by maximal nilpotence of $Q$ again. 
Let now $X$ be a normal subgroup of $\pi_1(Q)$ of finite index. Then 
$N_{G_1}(X)\ast \pi_2(Q)$ normalizes $X\ast \pi_2(Q)$ and as 
the latter is a normal subgroup of finite index in $Q$ one concludes that $X$ 
has finite index in $N_{G_1}(X)$. Hence $\pi_1(Q)$ is a Cartan subgroup of 
$G_1$.  Similarly, $\pi_2(Q)$ is a Cartan subgroup of $G_2$. 

Conversely, let $Q$ be a subgroup of $G$ of the form $Q_1 \ast Q_2$ 
for some Cartan subgroups $Q_i$ of $G_i$. Since each $Q_i$ is maximal nilpotent in $G_i$ 
it follows, considering projections as above, that $Q$ is maximal nilpotent in $G$. 
Let now $X$ be a normal subgroup of $Q$ of finite index. 
Then $\pi_i(N_G(X))$ normalizes the normal subgroup of finite index 
$\pi_i(X)$ of $Q_i$. Since $Q_i$ is a Cartan subgroup of $G_i$ it follows that $\pi_i(X)$ has finite 
index in $\pi_i(N_G(X))$. Finally, since $X\leq {\pi_1(X) \ast \pi_2(X)}\leq Q$, we get that 
$X$ has finite index in $N_G(X)$. 
\qed

The special case of a direct product in Lemma \ref{LemElemCartanCentProd} has also been observed 
in \cite[Chap. VI, \S4, Prop. 3]{MR0068552}. 

\begin{Corollary}\label{LemElemCartanDirProd}
Let $G=G_1\times \cdots \times G_n$ be a direct product of finitely many groups $G_i$. 
Then Cartan subgroups of $G$ are exactly of the form 
$Q_1\times \cdots \times Q_n$ where each $Q_i$ is a Cartan subgroup of $G_i$. 
\end{Corollary}

\section{Dimension and unions}\label{DimensionUnions}

In this section we work with a structure such that each nonempty definable set is equipped with a dimension 
in $\N$ satisfying the following axioms for any nonempty definable sets $A$ and $B$. 

\begin{itemize}
\item[(A1)]{\bf (Definability)} 
If $f$ is a definable function from $A$ to $B$, then the set 
$\{ b\in B~|~\dim(f^{-1}(b))=m \}$ is definable for every $m$ in $\N$. 
\item[(A2)]{\bf (Additivity)} 
If $f$ is a definable function from $A$ to $B$, whose fibers have constant 
dimension $m$ in $\N$, then $\dim(A)=\dim(\Image(f))+m$. 
\item[(A3)]{\bf (Finite sets)}
$A$ is finite if and only if $\dim(A)=0$.
\item[(A4)]{\bf (Monotonicity)}
$\dim(A\cup B)=\max(\dim(A),\dim(B))$.
\end{itemize}

In an o-minimal structure, definable sets are equipped with a finite dimension 
satisfying all these four axioms, by \cite[Chapter 4]{MR1633348} or \cite{MR961362}. 
Hence our reader only interested in groups definable in o-minimal structures may 
read all the following dimensional computations in the restricted context of such groups. 
But, as mentioned in the introduction, such computations are relevant in other 
contexts as well (groups of finite Morley rank, groups in supersimple theories of finite rank, 
groups definable over the p-adics...), and thus we will proceed with the mere axioms A1-4. 

Axioms A2 and A3 guarantee that if $f$ is a definable bijection between two definable sets $A$ and $B$, 
then $\dim(A)=\dim(B)$. Axiom A4 is a strong form of monotonicity in the sense that 
$\dim(A)\leq \dim(B)$ whenever $A\subseteq B$. 

\begin{Definition}
Let $\mathcal M$ be a first-order structure equipped with a dimension $\dim$ on 
definable sets and $X\subseteq Y$ two definable sets. We say that $X$ is: 
\begin{itemize}
\item[$(a)$]
{\em weakly generic} in $Y$ whenever $\dim(X)=\dim(Y)$. 
\item[$(b)$]
{\em generic} in $Y$ whenever $Y$ is a definable group covered by finitely many translates of $X$. 
\item[$(c)$]
{\em large} in $Y$ whenever $\dim(Y\setminus X)<\dim(Y)$. 
\end{itemize}
\end{Definition}

Clearly, genericity and largeness both imply weak genericity when the dimension satisfies 
axioms A1-4. 
If $G$ is a group definable in an o-minimal structure and $X$ is a large definable subset of $G$, then 
$X$ is generic: see \cite[Lemma 2.4]{MR961362} for a proof by compactness, and 
\cite[Section 5]{MR2668231} for a proof with precise bounds on the number of 
translates needed for genericity. 
In the sequel we are only going to use dimensional computations, hence the notions 
of weak genericity and of largeness. We are not going to use the notion of genericity 
(which is imported from the theory of stable groups in model theory), but we will make 
some apparently quite new remarks on genericity and Cartan subgroups in real Lie groups 
(Remark \ref{RemarkSL2Q1GenQ2NotGen} below). 

Our arguments for Cartan subgroups in groups definable in o-minimal structures will highly 
depend on computations of the dimension of their unions in the style of \cite{MR2225896}, 
and to compute the dimension of a union of definable sets we adopt the following 
geometric argument essentially due to Cherlin. 

Assume from now on that $X_a$ is a uniformly definable family of definable sets, 
with $a$ varying in a definable set $A$ and such that $X_a=X_{a'}$ if and only if $a=a'$. 
We have now a combinatorial geometry, 
where the set of points is $U:=\bigcup_{a\in A}X_a$, the set of lines is the set $\{X_a~|~a\in A\}$ 
in definable bijection with $A$, and the incidence relation is the natural one. 
The set of {\em flags} is then defined to be the subset of couples $(x,a)$ of 
$U\times A$ such that $x\in X_a$. By projecting the set of flags on the set of points, 
one sees with axiom A1 that for any $r$ such that $0\leq r \leq \dim(A)$, the set 
$$U_r:=\{x\in U~|~\dim(\{a\in A~|~x\in X_a\})=r\}$$
is definable. In particular, each subset of the form $[X_a]_r:=X_a\cap U_r$, i.e., the set of points $x$ 
of $X_a$ such the set of lines passing through $x$ has dimension $r$, is definable as well. 

\begin{Proposition}\label{GeneralRankComput}
In a structure equipped with a dimension satisfying axioms A1-2, let $X_a$ be a uniformly definable family 
of sets, with $a$ varing in a definable set $A$ and such that $X_a=X_{a'}$ if and only if $a=a'$. Suppose, 
for some $r$ such that $0\leq r\leq \dim(A)$, that $[X_a]_r$ is nonempty and that 
$\dim([X_a]_r)$ is constant as $a$ varies in $A$. Then 
$$\dim(U_r)+r=\dim(A)+\dim([X_a]_r).$$
\end{Proposition}
\proof
One can consider the definable subflag associated to $U_r=[\bigcup_{a\in A}X_a]_r$ in the 
point/line incidence geometry described above. 
By projecting this definable set on the set of points and on the set of lines respectively, 
one finds using axiom A2 of the dimension the desired equality as in 
\cite[\S2.3]{MR2225896}. 
\qed

Given a permutation group $(G,\Omega)$ and a subset $X$ of $\Omega$,
we denote by $N(X)$ and by $C(X)$ the {\em setwise} and
the {\em pointwise} stabilizer of $X$ respectively, that is
$G_{\{X\}}$ and $G_{(X)}$ in a usual
permutation group theory notation. We also denote by $X^{G}$ the set 
$\{x^{g}~|~(x,g)\in{X \times G}\}$, where $x^{g}$ denotes 
the image of $x$ under the action of $g$, as in the case of an action by conjugation. 
Subsets of the form $X^{g}$ for some $g$ in $G$ are also called
{\em $G$-conjugates} of $X$. Notice that the set $X^{G}$ can be seen, alternatively, as the 
union of $G$-orbits of elements of $X$, or also as the union of $G$-conjugates of $X$. 
When considering the action of a group on itself by conjugation,
as we will do below, all these terminologies and notations are the usual
ones, with $N(X)$ and $C(X)$ the {\em normalizer} and the {\em centralizer} of $X$ respectively. 

We shall now apply Proposition \ref{GeneralRankComput} in the context of permutation groups 
in a way much reminiscent of \cite[Fact 4]{MR2537672}. For that purpose we will need 
that the dimension is well defined on certain imaginaries, and for that purpose we will 
make the simplifying assumption that the theory considered eliminates such specific 
imaginaries. We recall that groups definable in o-minimal 
structures eliminate all imaginaries by \cite[Theorem 7.2]{MR2006422}, so these technical assumptions 
will always be verified in this context. (And our arguments are also valid in any context 
where the dimension is well defined and compatible in the relevant imaginaries.) 
For any quotient $X/{\sim}$ associated to an equivalence relation $\sim$ on a set $X$, 
we call {\em transversal} any subset of $X$ intersecting each equivalence class in exactly one point. 

\begin{Corollary}\label{FaitRankComput}
Let $(G,\Omega)$ be a definable permutation group in a structure 
equi\-pped with a dimension satisfying 
axioms A1-3, $X$ a definable subset of $\Omega$ such that $G/N(X)$ (right cosets) 
has a definable transversal $A$. 
Suppose that, for some $r$ between $0$ and $\dim(A)$, the definable subset 
$X_{r}:=\{x\in X~|~\dim(\{a\in A~|~x\in X^a\})=r\}$ is nonempty. Then 
$$\dim({X_{r}}^{G})=\dim(G)+\dim(X_{r})-\dim(N(X))-r.$$
\end{Corollary}
\proof
We can apply Proposition \ref{GeneralRankComput} with the uniformly definable family of $G$-conjugates 
of $X$, which is parametrized as $\{X^a~|~a\in A\}$ since $A$ is a definable transversal of $G/N(X)$. 
Notice that the sets $[X^a]_r$ are in definable bijection, as pairwise $G$-conjugates, and hence all 
have the same dimension. 
Notice also that $\dim(A)=\dim(G)-\dim(N(X))$ by the additivity of the dimension and its invariance under 
definable bijections. 
\qed

The following corollary, which is crucial in the sequel, can be compared to 
\cite[Corollary 5]{MR2537672}. 

\begin{Corollary}\label{CorHGenr=0}
Assume furthermore in Corollary \ref{FaitRankComput} that the dimension 
satisfies axiom A4, and that $\dim(G)=\dim(\Omega)$ and $\dim(X)\leq \dim(N(X))$. 
Then 
$$\dim(X^{G})=\dim(\Omega) {\mbox~if~and~only~if~} \dim(X_{0})=\dim(N(X))~(=\dim(X)).$$ 
In this case, ${X_0}^G$ is large in $X^G$. 
\end{Corollary}
\proof
If $\dim(X^{G})=\dim(\Omega)$, then one has for some $r$ as in
Corollary \ref{FaitRankComput} that $\dim({X_{r}}^{G})=\dim(\Omega)$ by axiom A4, and then
$$0\leq r =\dim(X_{r})-\dim(N(X))\leq \dim(X)-\dim(N(X))\leq 0$$
by monotonicity of the dimension, showing that all these quantities are equal to $0$. 
In particular $r=0$, and $\dim(X_{0})=\dim(N(X))$. 
Conversely, if $\dim(X_{0})=\dim(N(X))$,
then $\dim({X_{0}}^{G})=\dim(G)=\dim(\Omega)$ by Corollary \ref{FaitRankComput}. 

Assume now the equivalent conditions above are satisfied. 
The first part of the proof above shows that $\dim({X_r}^G)=\dim(X^G)$ 
($=\dim(\Omega)$) can occur 
only for $r=0$. Hence ${X_0}^G$ is large in $X^G$ by axiom A4 again. 
\qed

\begin{Remark}\label{RemarkBizarre}
In general it seems one cannot conclude also that $X_0$ is large in $X$ in 
Corollary \ref{CorHGenr=0}. One could imagine the (bizarre) configuration in which 
$\dim(X_r)=\dim(X)$ for some $r>0$; in this case $\dim({X_{r}}^{G})=\dim(\Omega)-r$. 
\end{Remark}

In the remainder we will always consider the action of a group $G$ on itself by conjugation, so the 
condition $\dim(G)=\dim(\Omega)$ will always be met in Corollary \ref{CorHGenr=0}. 
Then we can apply Corollary \ref{CorHGenr=0} with $X$ any normalizing coset of a 
definable subgroup $H$ of $G$, as commented in \cite[page 1064]{MR2537672}. 
More generally, we now see that we can apply it simultaneously to finitely many such cosets. 
We first elaborate on the notion of {\em generosity} defined in 
\cite{MR2225896} and \cite{MR2537672} in the finite Morley rank case. 

\begin{Definition}
\label{DefGenerous}
Let $X$ be a definable subset of a group $G$ definable in a structure equipped with a dimension 
satisfying axioms A1-4. We say that $X$ is 
\begin{itemize}
\item[$(a)$]
{\em weakly generous} in $G$ whenever $X^G$ is weakly generic in $G$. 
\item[$(b)$]
{\em generous} in $G$ whenever $X^G$ is generic in $G$. 
\item[$(c)$]
{\em largely generous} in $G$ whenever $X^G$ is large in $G$. 
\end{itemize}
\end{Definition}

\begin{Corollary}\label{CorHGenr=0WH}
Suppose $H$ is a definable subgroup of a group $G$ definable in a structure equipped with a dimension 
satisfying axioms A1-4, and suppose $W$ is a finite subset of $N(H)$ such that 
$G/N(WH)$ has a definable transversal. 
Then $WH$ is weakly generous in $G$ if and only if 
$$\dim([WH]_0)=\dim(N(WH)).$$ 
In this case, $[WH]_0^G$ is large in $[WH]^G$, and 
$\dim([WH]_0)=\dim(WH)=\dim(H)=\dim(N(WH))$. 
\end{Corollary}
\proof
Let $X=WH$. Since $W$ is finite, $X$ is definable. 
In order to apply Corollary \ref{CorHGenr=0}, one needs to check that 
$\dim(X)\leq \dim(N(X))$. Of course, the subgroup $H$ normalizes each coset $wH$, 
for each $w\in W\subseteq N(H)$, and in particular $H\leq N(WH)$. We get thus that 
$\dim(X)=\dim(WH)=\dim(H)\leq \dim(N(WH))=\dim(N(X))$. 

Now Corollary \ref{CorHGenr=0} gives our necessary and sufficient condition, and the largeness of 
$[WH]_0^G$ in $[WH]^G$. It also gives $\dim(X_0)=\dim(X)=\dim(N(X))$. We have seen already that 
$\dim(X)=\dim(H)$. 
\qed

The following lemma is a fundamental trick below. 

\begin{Lemma}\label{FondamentalLemma}
Let $G$ be a group definable in a structure equipped with a dimension satisfying 
axioms A1-4 and with the $dcc$. Let $X$ be a definable subset of $G$, $X_0$ 
the subset of elements of $X$ contained in only finitely many $G$-conjugates of $X$, 
and $U$ a definable subset of $X$ such that $U\cap X_0\neq \emptyset$. 
Then $N\o(U)\leq N(X)$. 
\end{Lemma}
\proof
As in \cite[Lemma 3.3]{MR2225896}, essentially via Fact \ref{ActionConOnFinite}$(a)$. 
\qed

\section{Cosets arguments}

Corollary \ref{CorHGenr=0WH} will be used at the end of this paper in certain arguments 
reminiscent of a theory of Weyl groups from \cite{MR2537672}. Since such specific arguments follow essentially 
from Corollary \ref{CorHGenr=0WH} we insert here, as a warm up, a short 
section devoted to them. 

\begin{Theorem}\label{TheoGenCosets}
Let $G$ be a group definable in a structure equipped with a dimension satisfying axioms A1-4 and with the $dcc$, 
$H$ a weakly generous definable subgroup of $G$, and $w$ an 
element normalizing $H$ and such that 
$G/N(H)$ has a definable transversal. 
Then one the following must occur: 
\begin{itemize}
\item[$(a)$]
The coset $wH$ is weakly generous in $G$, or 
\item[$(b)$]
The definable set $\{h^{w^{n-1}}h^{w^{n-2}}\cdots h~|~h\in H\}$ 
is not large in $H$ for any multiple $n$ of the (necessarily finite) order of $w$ modulo $H$. 
If $w$ centralizes $H$, then $\{h^n~|~h\in H\}$ is not large in $H$. 
\end{itemize}
\end{Theorem}
\proof
We proceed essentially as in \cite[Lemmas 11-12]{MR2537672}. 
Assume $wH$ not weakly generous. In particular $w\in {N(H)\setminus H}$ since 
$H$ is weakly generous by assumption. 
By Corollary \ref{CorHGenr=0WH}, $H_0$ is weakly generic in $N(H)$; in particular $H$ has finite 
index in $N(H)$. 
Of course, $N(wH)\leq N(H)$ since $H=\{ab^{-1}:a,b\in wH\}$, and one sees then that 
$N(wH)$ is exactly the preimage in $N(H)$ of the centralizer of $w$ modulo $H$. 
To summarize, $H \leq N(wH) \leq N(H)$, with $N(H)/H$ finite.  
In particular $w$ has finite order modulo $H$. Notice also at this stage that $G/N(wH)$ has 
a definable transversal (of the form $AX$ where $X$ is a definable transversal of $G/N(H)$ and 
$A$ is a definable transversal of the finite quotient $N(H)/N(wH)$). 
Since we assume $wH$ not weakly generous, 
Corollary \ref{CorHGenr=0WH} implies that 
$[wH]_0$ is not weakly generic in $wH$. In other words, the (definable) 
set of elements of the coset $wH$ contained 
in infinitely many $G$-conjugates of $wH$ is large in $wH$. 

Assume towards a contradiction $\{h^{w^{n-1}}h^{w^{n-2}}\cdots h~|~h\in H\}$ 
large in $H$ for $n$ a multiple of the finite order of $w$ modulo $H$. 
Let $\phi~:~wh\mapsto (wh)^{n}$ denote the definable map, 
from $wH$ to $H$, consisting of taking $n$-powers. As 
$$\phi(wH)=w^{n}\cdot \{h^{w^{n-1}}h^{w^{n-2}}\cdots h~|~h\in H\}$$ 
our contradictory assumption forces that $\phi(wH)$ must be large in $H$. 

Then $H_0\cap \phi(wH)$ must be weakly generic in $H$. 
Since the dimension can only get down when taking images by definable functions, 
$\phi^{-1}(H_0\cap \phi(wH))$ necessarily has to be weakly generic in the coset $wH$. 
Therefore one finds an element $x$ in the intersection of this preimage with the large 
subset $[wH]\setminus [wH]_0$ of elements of $wH$ 
contained in infinitely many $G$-conjugates of $wH$. 
Now since $w^n\in H$ and $N(wH)$ has finite index in $N(H)$ it follows that 
$\phi(x)=x^n$ belongs to infinitely many $G$-conjugates of $H$, a contradiction since 
$\phi(x)$ belongs to $H_0$. This proves our main statement in case $(b)$. 

For our last remark in case $(b)$, notice that when $w$ centralizes $H$ one has 
$\{h^{w^{n-1}}h^{w^{n-2}}\cdots h~|~h\in H\}=\{h^n~|~h\in H\}$. 
\qed

\begin{Corollary}\label{Cor1ThmGenCosets}
Suppose additionally in Theorem \ref{TheoGenCosets} that $w$ has order $n$ modulo $H$ and that $H$ 
is $n$-divisible ($n\geq 1$). Then one of the following must occur: 
\begin{itemize}
\item[$(a)$]
The coset $wH$ is weakly generous in $G$, or 
\item[$(b)$]
$C_H(w)$ is a proper subgroup of $H$. 
\end{itemize}
\end{Corollary}
\proof
Suppose that both alternatives fail. Then $\{h^n~|~h\in H\}$ is not large in $H$ by 
Theorem \ref{TheoGenCosets}, a contradiction since this set is $H$ by $n$-divisibility. 
\qed

The following corollary of Theorem \ref{TheoGenCosets} will be particularly adapted 
in the sequel to Cartan subgroups of groups definable in o-minimal structures. 

\begin{Corollary}\label{CorHGenCosetwHGen}
Suppose additionally in Theorem \ref{TheoGenCosets} that $H$ is definably connected and divisible 
and that $\<w\>H$ is nilpotent. Then the coset $wH$ is weakly generous in $G$. 
\end{Corollary}
\proof
This is clear if $w$ is in $H$, so we may assume $w\in {N(H)\setminus H}$. 
As above $w$ has finite order modulo $H=H\o$. 
By $dcc$ of the ambient group and \cite[Lemma 3.10]{BaroJalOteroCommutators}, 
the coset $wH$ contains a torsion element which 
commutes with $H=H\o$, and thus we may assume $C_H(w)=H$. By divisibility of $H=H\o$, 
$\{h^n~|~h\in H\}=H$ is large in $H$, and by Theorem \ref{TheoGenCosets} the coset $wH$ must 
be weakly generous in $G$. 
\qed

We will also use the following more specialized results in the same spirit, which apply as usual 
to nilpotent groups definable in o-minimal structures by \cite[Theorem 6.10]{MR2006422}. 

\begin{Lemma}\label{LemPhi(X)Gen}
Let $H$ be a nilpotent divisible group definable in a structure equipped with a dimension 
satisfying axioms A1-4, with the $dcc$, and with no infinite elementary abelian $p$-subgroups 
for any prime $p$. Let $\phi$ be the map consisting of 
taking $n$-th powers for some $n\geq 1$. 
If $X$ is a weakly generic definable subset of $H$, then $\phi(X)$ is weakly generic as well. 
\end{Lemma}
\proof
Considering the dimension, it suffices to show that $\phi$ has finite fibers. 
Suppose $a^n=b^n$ for some elements $a$ and $b$ in $H$. If $aZ(H)=bZ(H)$, 
then our assumption forces, with $a$ fixed, that 
$b$ can only vary in a finite set, as desired. Hence, working in $H/Z(H)$, it suffices to show 
that $a^n=b^n$ implies $a=b$. But by \cite[Lemma 3.10(a')]{BaroJalOteroCommutators} all 
definable sections of $H/Z(H)$ are torsion-free, and our claim follows easily by induction on the nilpotency 
class of $H/Z(H)$. 
\qed

\begin{Corollary}\label{CorLemPhi(X)Gen}
Let $Q$ be a nilpotent group definable in a structure equipped with a dimension satisfying 
axioms A1-4, with the $dcc$, and with no infinite elementary abelian $p$-subgroups 
for any prime $p$. Suppose $Q\o$ divisible, and let $a\in Q$, $n$ a multiple 
of the order of $a$ modulo $Q\o$, and $\phi$ the map consisting of taking $n$-th powers. 
If $X$ is a weakly generic definable subset of $aQ\o$, then $\phi(X)$ is a weakly generic subset of $Q\o$. 
\end{Corollary}
\proof
By \cite[Lemma 3.9]{BaroJalOteroCommutators}, we may assume that $a$ centralizes $Q\o$. 
Now for any $x\in Q\o$ we have $\phi(ax)=a^nx^n$. Hence, if $x$ varies in a weakly 
generic definable subset $X$ of $Q\o$, then $\phi(ax)$ also by Lemma \ref{LemPhi(X)Gen} in $H=Q\o$. 
\qed

\section{Generosity and lifting}

In the present section we study the behaviour of weak or large generosity when passing to 
quotients by definable normal subgroups. We continue with the mere axioms A1-4 of 
Section \ref{DimensionUnions} for the dimension, and with the existence of definable 
transversal for certain imaginaries to ensure that their dimensions is also well defined. 
As above, everything applies in particular to groups definable 
in o-minimal structures. 

\begin{Proposition}\label{PropGenLifting}
Let $G$ be a group definable in a structure equipped with a dimension satisfying axioms A1-4, 
$N$ a definable normal subgroup of $G$, $H$ a definable subgroup of $G$ containing $N$, 
and $Y$ a definable subset of $H$ large in $H$. 
Suppose also that $G/N$ and $G/N(H\setminus Y^H)$ have definable transversals. 
\begin{itemize}
\item[$(a)$]
If $H/N$ is weakly generous in $G/N$, then $Y$ is weakly generous in $G$. 
\item[$(b)$]
If $H/N$ is largely generous in $G/N$, then $Y$ is largely generous in $G$. 
\end{itemize}
\end{Proposition}
\proof
First note that $H^G$ is a union of cosets of $N$, since $N\leq H$ and $N\trianglelefteq G$. 
Hence the weak (resp. large) generosity  of $H/N$  in $G/N$ forces the weak (resp. large) 
generosity of $H$ in $G$. In any case, $\dim(H^G)=\dim(G)$. 

Replacing $Y$ by $Y^H$ if necessary, we may assume $H\leq N(Y)$ and $Y$ large in $H$. 

\begin{Claim}\label{ClaimLift2}
Let $Z=H\setminus Y$. Then $Z^G$ cannot be weakly generic in $H^G$. 
\end{Claim}
\proof
Suppose $Z^G$ weakly generic in $H^G$. Then $\dim(Z^G)=\dim(H^G)=\dim(G)$. 
Since $Z\subseteq H\subseteq N_G(Z)$, Corollary \ref{CorHGenr=0} yields 
$\dim(Z)=\dim(N_G(Z))$. In particular $\dim(Z)=\dim(H)$, a contradiction to the 
largeness of $Y$ in $H$.
\qed

$(a)$. Since $\dim(H^G)=\dim(G)$ and $H^G=Y^G\cup Z^G$, 
Claim \ref{ClaimLift2} yields $\dim(Y^G)=\dim(G)$.

$(b)$. In this case $H^G$ is large in $G$. Since 
$G=(G\setminus H^G) \sqcup (H^G\setminus Y^G) \sqcup Y^G$, Claim \ref{ClaimLift2} 
now forces $Y^G$ to be large in $G$.
\qed

\begin{Corollary}\label{CorLiftGenSubgroups}
Assume that $G$, $N$, $H$, and $Y$ are as in Proposition \ref{PropGenLifting}, and that 
$Y=Q^H$ for some largely generous definable subgroup $Q$ of $H$. 
\begin{itemize}
\item[$(a)$]
If $H/N$ is weakly generous in $G/N$, then so is $Q$ in $G$
\item[$(b)$]
If $H/N$ is largely generous in $G/N$, then so is $Q$ in $G$. 
\end{itemize}
\end{Corollary}
\proof
It suffices to apply Proposition \ref{PropGenLifting} with $Y=Q^H$, noticing that $Y^G=Q^G$. 
\qed

\begin{Corollary}\label{CorLiftCarterSubgroups}
Assume furthermore that $Q$ is a Carter subgroup of $H$ in Corollary \ref{CorLiftGenSubgroups}, and 
that $N_G(Q)/Q$ has a definable transversal. Then, in both cases $(a)$ and $(b)$, 
$Q$ is a Carter subgroup of $G$. 
\end{Corollary}
\proof
By definition, $Q$ is definable, definably connected, and nilpotent. So it suffices to check that 
$Q$ is a finite index subgroup of $N_G(Q)$. But in any case, 
it follows from the weak generosity of $Q$ in $G$ given in Corollary \ref{CorLiftGenSubgroups} 
and from Corollary \ref{CorHGenr=0WH} that $\dim(Q)=\dim(N_G(Q))$. 
Now axiom A3 applies. 
\qed

\section{Weakly generous nilpotent subgroups}\label{SectionWeakGenNilp}

In the present section we shall rework arguments from \cite{MR2225896} concerning weakly 
generous Carter subgroups. 
Throughout the section, $G$ is a group definable in a structure with a dimension satisfying 
axioms A1-4, and with the $dcc$. As in the preceding sections, 
everything applies in particular to groups definable in an o-minimal structure. 

\begin{Lemma}\label{FondamentalLemma_Nnilp}
Let $G$ be a group definable in a structure with a dimension satisfying 
axioms A1-4, and with the $dcc$. 
Let $H$ be a definable subgroup of $G$ such that $N\o(H)=H\o$, 
$H_0$ the set of elements of $H$ contained in 
only finitely many conjugates of $H$, and $N$ a definable nilpotent subgroup of $G$ such 
that $N\cap H_0$ is nonempty. Then $N\o\leq H\o$. 
\end{Lemma}
\proof
Let $U=N\cap H$. By assumption $U\cap H_0$ is nonempty, 
so by Lemma \ref{FondamentalLemma} 
$N\o(U)\leq N\o(H)=H\o$. In particular, $N\o_N(U)\leq {(N\cap H)\o}=U\o$, 
which shows that $U$ has finite index in $N_N(U)$. 
Now Lemma \ref{NormalizerCondition} shows that $U$ must have finite index in $N$, 
and in particular $U\o=N\o$. Hence, $N\o=(N\cap H)\o\leq H\o$. 
\qed

\begin{Corollary}\label{CorDescriptionQNilpwgen}
Let $G$ be a group definable in a structure with a dimension satisfying 
axioms A1-4, and with the $dcc$. 
Let $Q$ be a definable nilpotent weakly generous subgroup of $G$ such that 
$G/N(Q)$ has a definable transversal, and let $Q_0$ denote the set of elements of $Q$ contained 
in only finitely many conjugates of $Q$. Then: 
\begin{itemize}
\item[$(a)$]
For any definable nilpotent subgroup $N$ such that $N\cap Q_0\neq \emptyset$, we have 
$N\o\leq Q\o$. 
\item[$(b)$]
For any $g$ in $G$ such that $Q_0\cap{Q}^g \neq \emptyset$, we have that 
$Q\o={[Q\o]}^g$.  
\end{itemize}
\end{Corollary}
\proof
$(a)$. As $Q$ is weakly generous, we have $N\o(Q)=Q\o$ by Corollary \ref{CorHGenr=0WH}. 
Hence Lemma \ref{FondamentalLemma_Nnilp} gives $N\o\leq Q\o$. 
$(b)$. Item $(a)$ applied with $N=Q^g$ yields $[Q\o]^g=[Q^g]\o\leq Q\o$. Now applying 
Lemma \ref{NormalizerCondition} shows that $[Q\o]^g$ cannot be of infinite index in $Q\o$ 
(as otherwise we would contradict that $N\o(Q)=Q\o$), and thus $[Q\o]^g=Q\o$. 
\qed

\begin{Corollary}\label{CorDescriptionQwgenCarter}
Suppose in addition in Corollary \ref{CorDescriptionQNilpwgen} that $Q$ is a 
Carter subgroup of $G$. Then, for any $g\in Q_0$ and any definably 
connected definable nilpotent subgroup $N$ containing $g$, we have $N\leq Q$. 
In particular, $Q$ is the unique maximal definably connected definable nilpotent 
subgroup containing $g$, and the distinct conjugates of $Q_0$ are indeed disjoint, 
forming thus a partition of a weakly generic subset of $G$. 
\end{Corollary}
\proof
It suffices to apply Corollary \ref{CorDescriptionQNilpwgen}. 
\qed

As a result one also obtains the following general theorem, 
which can be compared to the main result of \cite{MR2225896}. 

\begin{Theorem}\label{TheoConjLargelyGenCarter}
Let $G$ be a group definable in a structure with a dimension satisfying 
axioms A1-4, and with the $dcc$. 
Then $G$ has at most one conjugacy class of largely generous Carter subgroups $Q$ 
such that $G/N(Q)$ has a definable transversal. 
If such a Carter subgroup exists, then the set of elements contained in a unique conjugate of that 
Carter subgroup is large in $G$. 
\end{Theorem}
\proof
Let $P$ and $Q$ be two largely generous Carter subgroups of $G$. 
We want to show that $P$ and $Q$ are conjugate. 
We have ${P_0}^G$ and ${Q_0}^G$ large in $G$ by Corollary \ref{CorHGenr=0WH}. 
Since the intersection of two large sets is nontrivial (and in fact large as well), 
we get that ${P_0}^G \cap {Q_0}^G$ is nonempty, so after conjugation we may thus assume 
$P_0 \cap Q_0$ nonempty. But then Corollary \ref{CorDescriptionQwgenCarter} gives $P=Q$. 

Our last claim follows also from Corollary \ref{CorDescriptionQwgenCarter}. 
\qed

\section{On groups definable in o-minimal structures}\label{SectionFactsGpsInoMinStruct}

We shall now collect results specific to groups definable in o-minimal structures which are 
needed in the sequel. 
We recall that groups definable in o-minimal structures satisfy the $dcc$ on definable subgroups 
\cite[Remark 2.13]{MR961362}, and o-minimal structures are equipped with a dimension 
satisfying axioms A1-4 considered in the previous sections 
\cite[Chapter 4]{MR1633348}. As commented before, we can freely apply all the results of the 
preceding sections to the specific case of groups definable in an o-minimal structure. 
We also recall that all the 
technical assumptions on the existence of transversals in 
Sections \ref{DimensionUnions}-\ref{SectionWeakGenNilp} are satisfied, since groups definable 
in o-minimal structures eliminate all imaginaries by \cite[Theorem 7.2]{MR2006422}. 
As mentioned already in the introduction, we consider only groups $G$ {\em definable} in an o-minimal 
structure, but \cite[Theorem 7.2]{MR2006422} also allows one to consider any group of the form 
$K/L$, where $L\trianglelefteq K \leq G$ are definable subgroups, as definable. 

\begin{Fact}\label{FactCommutators}
{\bf \cite[\S6]{BaroJalOteroCommutators}}
Let $G$ be a group definable in an o-minimal structure, with $G\o$ solvable, and $A$ and $B$ two 
definable subgroups of $G$ normalizing each other. Then $[A,B]$ is definable, and definably connected whenever $A$ and $B$ are. 
\end{Fact}

Any group $G$ definable in an o-minimal structure has a largest normal nilpotent subgroup 
$F(G)$, which is also definable \cite[Fact 3.5]{BaroJalOteroCommutators}, 
and a largest normal solvable subgroup $R(G)$, which is also definable 
\cite[Lemma 4.5]{BaroJalOteroCommutators}. 

\begin{Fact}\label{FactGenGDefConRes}
Let $G$ be a definably connected solvable group definable in an o-minimal structure. 
\begin{itemize}
\item[$(a)$] 
\cite[Theorem 6.9]{MR2006422} $G'$ is nilpotent. 
\item[$(b)$]
\cite[Proposition 5.5]{BaroJalOteroCommutators} $G'\leq F\o(G)$. In particular $G/F\o(G)$ and $G/F(G)$ are divisible abelian groups. 
\item[$(c)$]
\cite[Corollary 5.6]{BaroJalOteroCommutators} If $G$ is nontrivial, then $F\o(G)$ is nontrivial. 
In particular $G$ has an infinite abelian characteristic definable subgroup. 
\item[$(d)$]
\cite[Lemma 3.6]{BaroJalOteroCommutators} If $G$ is nilpotent and $H$ is an infinite normal subgroup of $G$, then $H\cap Z(G)$ is infinite. 
\end{itemize}
\end{Fact}

If $H$ and $G$ are two subgroups of a group with $G$ normalizing $H$, then a 
{\em $G$-minimal} subgroup of $H$ is an infinite $G$-invariant definable subgroup of $H$, 
which is minimal with respect to these properties (and where definability refers to the fixed 
underlying structure, as usual). If $H$ is definable and satisfies the $dcc$ on definable subgroups, then 
$G$-minimal subgroups of $H$ always exist. As the definably 
connected component of a definable subgroup is a definably characteristic subgroup, 
we get also in this case that any $G$-minimal subgroup of $H$ should be definably connected. 

\begin{Lemma}\label{LemGMinimalSubgGpRes}
Let $G$ be a definably connected solvable group definable in an o-minimal structure, 
and $A$ a $G$-minimal subgroup of $G$. Then $A\leq Z\o(F(G))$, and 
$C_G(a)=C_G(A)$ for every nontrivial element $a$ in $A$. 
\end{Lemma}
\proof
By Fact \ref{FactGenGDefConRes}$(c)$, $A$ has an infinite characteristic 
abelian definable subgroup. 
Therefore the $G$-minimality of $A$ forces $A$ to be abelian. 
In particular, $A\leq F(G)$. Since $A$ is normal in $F(G)$, Fact \ref{FactGenGDefConRes}$(d)$
and the $G$-minimality of $A$ now force that $A\leq Z(F(G))$. Since $A$ is definably connected, 
we have indeed $A\leq Z\o(F(G))$. 

Now $F(G)\leq C_G(A)$, and $G/C_G(A)$ is definably isomorphic to a quotient of $G/F(G)$. 
In particular $G/C_G(A)$ is abelian by Fact \ref{FactGenGDefConRes}$(b)$. 
If $A\leq Z(G)$, then clearly $C_G(a)=C_G(A)$ ($=G$) for every $a$ in $A$, and thus we may assume 
$G/C_G(A)$ infinite. Consider the semidirect product $A\rtimes (G/C_G(A))$. 
Since $A$ is $G$-minimal, $A$ is also $G/C_G(A)$-minimal. 
Now an o-minimal version of Zilber's Field Interpretation Theorem 
for groups of finite Morley rank \cite[Theorem 2.6]{MR1779482} applies directly to 
$A\rtimes (G/C_G(A))$. It says that there is an infinite 
interpretable field $K$, with $A\simeq K_+$ and $G/C_G(A)$ an infinite subgroup of $K^{\times}$, and 
such that the action of $G/C_G(A)$ on $A$ corresponds to scalar multiplication. In particular, 
$G/C_G(A)$ acts {\em freely} (or {\em semiregularly} in another commonly used terminology) 
on $A\setminus \{1\}$. This means exactly 
that for any nontrivial element $a$ in $A$, $C_G(a)\leq C_G(A)$, i.e., $C_G(a)=C_G(A)$. 
\qed

For definably connected groups definable in an o-minimal structure which are not 
solvable, our study of Cartan subgroups will make heavy use of the main theorem about 
groups definable in o-minimal structures. It can be summarized 
as follows, compiling several papers to which we will refer immediately after the statement. 
Recall that a group is {\em definably simple} if the only definable normal subgroups are the 
trivial and the full subgroup. 

\begin{Fact}\label{TheoCherlinZilberConjoMin}
Let $G$ be a definably connected group definable in an o-minimal structure $\mathcal{M}$. 
Then 
$$G/R(G)=G_1\times \cdots \times G_n$$ 
where each $G_i$ is a definably simple infinite definable group. 
Furthermore, for each $i$, there is an $\mathcal{M}$-definable real closed field $R_i$ such that 
$G_i$ is $\mathcal{M}$-definably isomorphic to a semialgebraically connected semialgebraically simple 
linear semialgebraic group, definable in $R_i$ over the subfield of real algebraic numbers of $R_i$. 

Besides, for each $i$, either 
\begin{itemize}
\item[$(a)$] 
$\langle G_i, \cdot \rangle$ and 
$\langle R_i(\sqrt{-1}), +, \cdot \rangle$ are bi-interpretable; in this case $G_i$ is definably
isomorphic in $\langle G_i, \cdot \rangle$ to the $R_i(\sqrt{-1})$-rational points of a linear algebraic 
group, or
\item[$(b)$] 
$\langle G_i, \cdot \rangle$ and $\langle R_i, +, \cdot\rangle$ are 
bi-interpretable; in this case $G_i$ is definably isomorphic in 
$\langle G_i, \cdot \rangle$ to the connected component of the $R_i$-rational points of an algebraic 
group without nontrivial normal algebraic subgroups defined over $R_i$.
\end{itemize}
\end{Fact}

The description of $G/R(G)$ as direct product
of definably simple definable groups can be found in 
\cite[4.1]{MR1707202}. 
The second statement about definably simple groups is in \cite[4.1 \& 4.4]{MR1707202}, with the remark 
concerning the parameters in the proof of \cite[5.1]{MR1873380}. 
The final alternative for each factor, essentially between the complex case and the real case, 
is in \cite[1.1]{MR1779482}. 

When applying Fact \ref{TheoCherlinZilberConjoMin} in the sequel we will also use the following. 

\begin{Remark}\label{DimM=DimR}
Let $\mathcal{M}$ be an o-minimal structure, $R$ a real closed field definable in $\mathcal{M}$, and $X$ an 
$R$-definable subset of some $R^n$. Then $\dim_{\mathcal{M}}(X)=\dim_R(X)$. 
\end{Remark}
\proof
By o-minimality, $\mathcal{M}$ is a geometric structure 
\cite[Definition 3.2]{MR1779482}. Moreover, since $\dim_{\mathcal{M}}(R)=1$ 
by \cite[Proposition 3.11]{MR961362} we deduce that $R$ itself is $R$-minimal in 
the sense of \cite[Definition 3.3]{MR1779482}. Hence by \cite[Lemma 3.5]{MR1779482} 
we have that $\dim_{\mathcal{M}}(X)=\dim_{\mathcal{M}}(R)\dim_{R}(X)=\dim_{R}(X)$. 
\qed

We finish the present section with specific results about definably compact groups 
which might be used when such specific groups are involved in the sequel. 

\begin{Fact}\label{FactCompactGps}
Let $G$ be a definably compact definably connected group definable in an o-minimal structure. 
\begin{itemize}
\item[$(a)$] 
\cite[Corollary 5.4]{MR1729742} Either $G$ is abelian or $G/Z(G)$ is semisimple. 
In particular, if $G$ is solvable, then it is abelian. 
\item[$(b)$] 
\cite[Proposition 1.2]{MR2178762} 
$G$ is covered by a single conjugacy class of a definably connected definable
abelian subgroup $T$ such that $\dim(T)$ is maximal among dimensions of abelian definable 
subgroups of $G$.
\end{itemize}
\end{Fact}

For a variation on Fact \ref{FactCompactGps}$(b)$, see also \cite[Corollary 6.13]{MR2401858}. 
With Fact \ref{FactCompactGps} we can entirely clarify properties of Cartan subgroups 
in the specific case of definably compact groups definable in o-minimal structures, 
with a picture entirely similar to that in compact real Lie groups. 

\begin{Corollary}\label{CorStructCompactGps}
Let $G$ be a definably compact definably connected group definable in an o-minimal structure. Then 
Cartan subgroups $T$ of $G$ exist and are abelian, definable, definably connected, and conjugate, and 
$G=T^G$. 
\end{Corollary}
\proof
Let $T$ be a definably connected abelian subgroup as in Fact \ref{FactCompactGps}$(b)$. 
Since $G=T^G$, $T$ is in particular weakly generous, and 
thus of finite index in its normalizer by Corollary \ref{CorHGenr=0WH}. 
Hence $T$ is a Carter subgroup of $G$. 
Since $G=T^G$ again, and $t\in T \leq C\o(t)$ for every $t\in T$, we have the property that 
$g\in C\o(g)$ for every $g$ in $G$. 

We now prove our statement by induction on $\dim(G)$. 
By Lemma \ref{Carter-Cartan-Gen}, $T\leq Q$ for some Cartan subgroup such that 
$Q\o=T$. This takes care of the existence of Cartan subgroups of $G$, and their definability 
follows from Lemma \ref{Carter-Cartan-Gen}$(a)$. We also have $G=T^G$. 
We now claim that $T=Q$. Otherwise, $T=Q\o<Q$, and we find by Fact \ref{StructNilpGps} an element 
$a$ in $Q\setminus T$ centralizing $T$. Since $a\in T^g$ for some $g\in G$, we have 
$T$ and $T^g$ in $C\o(a)$. Now the Carter subgroups $T$ and $T^g$ of $C\o(a)$ are 
conjugate by an element of $C\o(a)$, obviously if $C\o(a)=G$ and by induction otherwise. 
Since $a\in T^g \leq C\o(a)$, we get $a\in T$, a contradiction. 
Hence $T=Q$ is a Cartan subgroup of $G$. 

It remains just to show that 
Cartan subgroups of $G$ are conjugate. 
Let $Q_1$ be an arbitrary Cartan subgroup of $G$, and $z$ a nontrivial 
element of $Z(Q_1)$ (Fact \ref{StructNilpGps} and 
Fact \ref{FactGenGDefConRes}$(d)$). We also have $z\in T^g$ for some $g\in G$, and thus 
$Q_1$, $T^g\leq C(z)$. If $C\o(z)<G$, the induction hypothesis applied in $C\o(z)$ yields the 
conjugacy of $Q_1\o$ and of $T$, giving also $Q_1=Q_1\o$ by maximal nilpotence of $T$. 
So we may assume $z\in Z(G)$. If $Z(G)$ is finite, then $G/Z(G)$ has a trivial center by 
Fact \ref{ActionConOnFinite}$(b)$, and the previous argument applied in $G/Z(G)$, together 
with Lemma \ref{LemElemCartanModZn}$(a)$, yields the conjugacy of $Q_1$ and $T$. 
Remains the case $Z(G)$ infinite: then applying the induction hypothesis in 
$G/Z(G)$, and using Lemma \ref{LemElemCartanModZn}$(a)$, also gives the conjugacy of 
$Q_1$ and $T$. This completes our proof. 
\qed

We have seen in the proof of Corollary \ref{CorStructCompactGps} that the 
``maximal definable-tori" $T$ of Fact \ref{FactCompactGps}$(b)$ must be Cartan subgroups, 
and then the two types of subgroups coincide by the conjugacy of Cartan subgroups. 
We note that the conjugacy of 
the ``maximal definable-tori" $T$ as in Fact \ref{FactCompactGps}$(b)$ was also shown 
in \cite{MR2178762}. 
Besides, we note that the maximal nilpotence of a Cartan subgroup $T$ of a group $G$ 
always implies that $C_G(T)=Z(T)$. In particular, in Corollary \ref{CorStructCompactGps}, 
$C(T)=T$ and the ``Weyl group" $W(G,T):=N(T)/C(T)$ acts faithfully on $T$. 

Finally, we take this opportunity to mention, parenthetically, a refinement of Fact \ref{FactCompactGps}$(a)$. 

\begin{Corollary}
Let $G$ be a definably compact definably connected group definable in an o-minimal structure. 
Then $R(G)=Z(G)$. 
\end{Corollary}
\proof
By Fact \ref{FactCompactGps}$(a)$ and \cite[Lemma 3.13]{BaroJalOteroCommutators}. 
\qed

\section{The definably connected solvable case}\label{SectionSolvableCase}

In the present section we are going to prove the following theorem. 

\begin{Theorem}\label{MainTheoSolvGps}
Let $G$ be a definably connected solvable group definable in an o-minimal structure. 
Then Cartan subgroups of $G$ exist and are conjugate, and they are definably connected and selfnormalizing. 
Moreover, they are largely generous in the following strong sense: for any Cartan subgroup $Q$, the 
(definable) set of elements of $Q$ contained in a unique conjugate of $Q$ is large in $Q$ and largely 
generous in $G$.  
\end{Theorem}

We first look at the minimal configuration for our analysis which can be 
thought of an abstract analysis of Borel subgroups of $\SL_2$ (over $\C$ or $\R$), first studied by Nesin 
in the case of groups of finite Morley rank \cite[Lemma 9.14]{MR1321141}. 

\begin{Lemma}\label{ConfigMinimaleRes}
Let $G$ be a definably connected solvable group definable in an o-minimal structure, 
with $G'$ a $G$-minimal subgroup and $Z(G)$ finite. Then $G=G'\rtimes Q$ for some (abelian) selfnormalizing 
definably connected definable largely generous complement $Q$, and 
any two complements of $G'$ are $G'$-conjugate. More precisely, we also have: 
\begin{itemize}
\item[$(a)$]
$F(G)={Z(G)\times G'}=C_G(G')$. 
\item[$(b)$]
For any $x$ in $G\setminus F(G)$, $xG'=x^{G'}$, $G=G'\rtimes C(x)$, and $C(x)$ is the unique 
conjugate of $C(x)$ containing $x$. 
\end{itemize}
\end{Lemma}
\proof
We elaborate on the proof given in \cite[Theorem 3.14]{MR2436138} in the finite Morley rank case. 
Since $Z(G)$ is finite, the definably connected group $G$ is not nilpotent by 
Fact \ref{FactGenGDefConRes}$(d)$, and in particular $C_G(G')<G$. By $G$-minimality of $G'$ and 
Lemma \ref{LemGMinimalSubgGpRes}, $G'\leq Z\o(F(G))$ and $C_G(a)=C_G(G')$ 
for every non-trivial element $a$ of $G'$. 

For any element $x$ in $G\setminus C_G(G')$, we now show that $Q:=C_G(x)$ is a required 
complement of $G'$. 
Since $x\notin C_G(G')$, $C_{G'}(x)=1$ and in particular $\dim(x^G)\geq \dim(G')$. 
On the other hand, $x^G\subseteq xG'$ as $G/G'$ is abelian, and it follows that $\dim(x^G)=\dim(G')$, or in other 
words that $\dim(G/Q)=\dim(G')$. Since $Q\cap G'=1$, the definable subgroup 
$G'\rtimes Q$ has maximal dimension in $G$, and since $G$ is definably connected we get that 
$G=G'\rtimes Q$. Of course $Q\simeq G/G'$ is abelian, and definably connected as $G$ is. 
We also see that $N_{G'}(Q)=C_{G'}(Q)=1$, since $C_{G'}(x)=1$, and thus the definable subgroup 
$Q=C_G(x)$ is selfnormalizing. 

$(a)$. The finite center $Z(G)$ is necessarily in $Q=C_G(x)$ in the previous paragraph, and in particular 
${Z(G)\cap G'}=1$. Since $G=G'\rtimes Q$ and $Q$ is abelian, $C_{Q}(G')\leq Z(G)$, 
and since $G'\leq Z(F(G))$ one gets $Z(G)\times G'\leq F(G)\leq C_G(G')\leq Z(G)\times G'$, proving 
item $(a)$. 

$(b)$. Let again $x$ be any element in $G\setminus F(G)$. The map 
$G'\to G'\colon u\mapsto [x,u]$ is a definable group homomorphism since $G'$ is abelian, with 
trivial kernel as $C_{G'}(x)=1$, and an isomorphism onto $G'$ since the latter is definably connected. 
It follows that any element of the form $xu'$, for $u'\in G'$, has the form 
$xu'=x[x,u]=x^u$ for some $u\in G'$, i.e., $xG'=x^{G'}$. 

Next, notice that any complement $Q_1$ of $G'$ is of the form $Q_1=C_G(x_1)$ for any 
$x_1\in Q_1\setminus Z(G)$. Indeed, $x_1\not\in Z(G)$ and $Q_1$ abelian imply 
$x_1\not\in C_G(G')$, and as above $C_G(x_1)$ is a definably connected complement of $G'$ 
containing $Q_1$, and comparing the dimensions we get $Q_1=C_G(x_1)$.

Moreover, if $Q_1=C_G(x_1)$  and $Q_2=C_G(x_2)$ are two complements of $G'$, 
we can always choose $x_1$ and $x_2$ in the same $G'$-coset; then they are $G'$-conjugate, 
as well as $Q_1$ and $Q_2$.
It is also now clear that, for any $x\in {G\setminus F(G)}$, $C_G(x)$ is the unique complement of $G'$ containing $x$, proving item $(b)$. 

It is clear from item $(b)$ that two complements of $G'$ are $G'$-conjugate, and that such complements 
are largely generous in $G$. 
\qed

\begin{Corollary}\label{CorConfigMinimaleRes} 
Let $G$ be a group as in Lemma \ref{ConfigMinimaleRes}. Then: 
\begin{itemize}
\item[$(a)$]
If $X$ is an infinite subgroup of a complement $Q$ of $G'$, then $N_G(X)=Q$ and 
$N_G(X)\cap G'=1$.
\item[$(b)$]
If $X$ is a nilpotent subgroup of $G$ not contained in $F(G)$, then $X$ is in an abelian 
complement of $G'$. 
\item[$(c)$]
Complements of $G'$ in $G$ are both Carter and Cartan subgroups of $G$, 
and all are of this form. 
\end{itemize} 
\end{Corollary}
\proof
$(a)$. We have $Q\leq N_G(X)$, and thus $N_G(X)=N_{G'}(X)\rtimes Q$. But 
$[N_{G'}(X),X]\leq {N_{G'}(X)\cap X}=1$ since $Q\cap G'=1$. In view of Lemma \ref{ConfigMinimaleRes}, 
and since $X$ is infinite, the only possibility is that $N_{G'}(X)=1$. Hence $N_G(X)=Q$, which is disjoint 
from $G'$. 

$(b)$. $X$ contains an element $x$ outside of $F(G)=C_{G}(G')$. 
Replacing $X$ by its definable hull $H(X)$ and using Fact \ref{FactGenGdcc}$(b)$, 
we may assume without loss that $X$ is definable. As in the proof of 
Lemma \ref{ConfigMinimaleRes}, ${X\cap G'}=\{[x,u]~|~u\in {X\cap G'}\}$, and the nilpotency of $X$ forces that 
${X\cap G'}=1$. Hence $X$ is abelian, and in the complement $C(x)$ of $G'$. 

$(c)$. Complements of $G'$ are selfnormalizing Carter subgroups by 
Lemma \ref{ConfigMinimaleRes}, and thus also Cartan subgroups by Lemma \ref{Carter-Cartan-Gen}. 
Conversely, one sees easily that a Carter or a Cartan subgroup of $G$ cannot be contained 
in $F(G)$, and then must be a complement of $G'$ by item $(b)$. 
\qed

Crucial in our proof of Theorem \ref{MainTheoSolvGps}, the next point shows that any 
definably connected nonnilpotent solvable group has a quotient 
as in Lemma \ref{ConfigMinimaleRes}. 

\begin{Fact}\label{GModWSL2Type}
{\bf (Cf. \cite[Proposition 3.5]{MR1765797})}
Let $G$ be a definably connected nonnilpotent solvable group definable in an o-minimal structure. 
Then $G$ has a definably connected definable normal subgroup $N$ such that $(G/N)'$ is $G/N$-minimal 
and $Z(G/N)$ is finite. 
\end{Fact}
\proof
The proof works formally exactly as in \cite[Proposition 3.5]{MR1765797} in the finite Morley rank case. 
All facts used there about groups of finite Morley rank have their formal analogs in 
Fact \ref{FactGenGDefConRes}$(a)$ and Lemma \ref{LemGMinimalSubgGpRes} in the o-minimal case. 
We also use the fact that lower central series and derived series of definably connected solvable groups 
definable in o-minimal structures are definable and definably connected, which follows from 
Fact \ref{FactCommutators} here. 
\qed

We now pass to the proof of the general Theorem \ref{MainTheoSolvGps}. 
At this stage we could follow the analysis by {\em abnormal} subgroups 
of \cite{Carter61/62} in finite solvable groups, 
as developed in the case of infinite solvable groups of finite Morley rank in \cite{MR1765797}. 
However we provide a more conceptual proof of Theorem \ref{MainTheoSolvGps}, mixing the use of 
Fact \ref{GModWSL2Type} with our general genericity arguments, in particular of 
Section \ref{SectionWeakGenNilp}. 
We note that the proof of Theorem \ref{MainTheoSolvGps} 
we give here would work equally in the finite Morley rank case (in that case there is no elimination 
of imaginaries but the dimension is well defined on imaginaries), providing a somewhat more 
conceptual proof of the analogous theorem in \cite{MR1765797} in that case. 

\medskip
\noindent
{\bf Proof of Theorem \ref{MainTheoSolvGps}.} 
We proceed by induction on $\dim(G)$. Clearly a minimal counterexample $G$ has to be nonnilpotent, and then 
has a definably connected definable normal subgroup $N$ as in Fact \ref{GModWSL2Type}. 
In what follows we use the notation ``$\overline{\phantom{H}}$" to denote quotients by $N$. 
Notice that $\overline{G}$ is necessarily infinite in Fact \ref{GModWSL2Type}, and $N$ is a subgroup of infinite index 
in $G$. 

\begin{Claim}\label{ClaimTheoRes1}
$G$ contains a definably connected and selfnormalizing Cartan subgroup $Q$ 
which is largely generous in the following sense: 
the (definable) set of elements of $Q$ contained in a unique conjugate of $Q$ is large in $Q$ and largely generous 
in $G$. 
\end{Claim}
\proof
Let $H$ be a definable subgroup of $G$ containing $N$ such that $\overline{H}$ is a selfnormalizing 
largely generous Carter subgroup of $\overline{G}$ as in Lemma \ref{ConfigMinimaleRes}. 
Notice that $H$ is definably connected since $\overline{H}$ and $N$ are. 
As $\overline{G}'$ is infinite, $\dim(\overline{H})<\dim(\overline{G})$, and $\dim(H)<\dim(G)$. 
We can thus apply the induction hypothesis in $H$, 
and assume that $H$ contains a definably connected and 
selfnormalizing Cartan subgroup $Q$ with the strong large generosity property: 
the set of elements of $Q$ contained in a unique $H$-conjugate of $Q$ is large in $Q$ 
and largely generous in $H$. We will show that $Q$ is the required subgroup. 

First note that $Q$, being definably connected, is a largely generous Carter subgroup of $H$. By 
Corollaries \ref{CorLiftGenSubgroups} and \ref{CorLiftCarterSubgroups}, $Q$ must be a largely 
generous Carter subgroup of $G$. 
We now show that $Q$ is selfnormalizing in $G$. Notice that $Q$ has an infinite image in 
$\overline{H}$, since it is largely generous in $H$ and $N$ is normal and proper in $H$. 
If $x\in N_G(Q)$, then $\overline{x}\in N_{\overline{G}}(\overline{Q})=\overline{H}$ by 
Corollary \ref{CorConfigMinimaleRes}$(a)$, and since $Q$ is selfnormalizing in $H$ it follows that 
$x\in N_H(Q)=Q$. Thus $Q$ is selfnormalizing in $G$. By Lemma \ref{Carter-Cartan-Gen}, 
$Q$ is also a Cartan subgroup of $G$. 

It remains just to show the largeness issue. 
Let $Q_{0}$ denote the set of elements of $Q$ contained in a unique $H$-conjugate of $Q$. 
We know that $Q_0$ is large in $Q$ and that $[Q_0]^H$ is large in $H$, 
so $[Q_0]^G$ ($=[[Q_0]^H]^G$) is large in $G$ by Proposition \ref{PropGenLifting}. This shows that 
$Q$ is largely generous in $G$, and thus it remains only to show it is in the strong sense of our claim. 
For that purpose, one easily sees that it is enough to show that the subset $X$ of elements of $Q_0$ 
contained in a unique {\em $G$-conjugate} of $Q$ is still large in $Q_0$, given the large partition of $G$ 
as in Corollary \ref{CorDescriptionQwgenCarter} and Theorem \ref{TheoConjLargelyGenCarter} 
(see also Proposition \ref{GeneralRankComput}). Since $Q$ is largely generous in $H$ and 
the preimage $L$ in $H$ of $F(\overline{G})$ is normal and proper in $H$, we get that $Q\nleq L$, 
and thus it suffices to show that $Q_0\setminus X$ is in $L$. 
Suppose towards a contradiction that an element $x$ in $Q_0$ and not in $L$ is in $Q^g$ for some 
$g$ not in $N_G(Q)$. Looking at images in $\overline{G}$ and since 
$\overline{x}\in{\overline{H}\setminus Z(\overline{G})}$, we then see with 
Lemma \ref{ConfigMinimaleRes} that $\overline{g}\in{N_{\overline{G}}(\overline{H})}={\overline{H}}$, 
and thus $g\in H$. 
Then $x\in Q\cap Q^g$ for some $g\in {H\setminus N_H(Q)}$, a contradiction since $x$ is in a unique 
$H$-conjugate of $Q$. 
This completes our proof of Claim \ref{ClaimTheoRes1}. 
\qed

\begin{Claim}\label{ClaimTheoRes2}
Carter subgroups of $G$ are conjugate. 
\end{Claim}
\proof
There are indeed at this stage two quick ways to argue for the conjugacy of Carter subgroups, 
either by quotienting by a $G$-minimal subgroup of $G$ as in 
\cite[Proof of Theorem 3.11]{MR2436138}, 
or still looking at the quotient $\overline{G}$. Since we have already 
used $\overline{G}$ for the existence of 
a largely generous Carter subgroup we keep on this second line of arguments. 

Let $Q_1$ be an arbitrary Carter subgroup of $G$.  By Theorem \ref{TheoConjLargelyGenCarter}, 
it suffices to prove that $Q_1$ is a largely generous Carter subgroup of $G$. 
Let $L$ be the preimage of $[\overline{G}]'$ in $G$; notice that $L$ is definably 
connected as $[\overline{G}]'$ and $N$ are. If $Q_1 \leq L$, then a Frattini Argument applied in $L$, 
using the induction hypothesis in $L$, gives 
$G=L\cdot N_G(Q_1)$, and since $Q_1$ is a Carter subgroup this gives that 
$L$ has finite index in $G$, a contradiction. 
Therefore $Q_1\nleq L$, and since $\overline{Q_1}$ is definably connected 
we also get $\overline{Q_1}\nleq F(\overline{G})$ by Lemma \ref{ConfigMinimaleRes}$(a)$. 
In particular, by Corollary \ref{CorConfigMinimaleRes}$(b)$, $Q_1$ is contained in a 
definably connected definable subgroup $H$ as in the proof of Claim \ref{ClaimTheoRes1}. 
Since $H<G$, the induction hypothesis 
applies in $H$, and thus $Q_1$ must be conjugate in $H$ to a largely generous 
Carter subgroup $Q$ of $H$.  
In particular, by the proof of Claim \ref{ClaimTheoRes1}, $Q_1$ is a largely generous Carter 
subgroup of $G$, as required.
\qed

The Cartan subgroup $Q$ provided by Claim \ref{ClaimTheoRes1} is also a Carter subgroup 
by definable connectedness and Lemma \ref{Carter-Cartan-Gen}$(a')$. 
If $Q_1$ is an arbitrary Cartan subgroup, then $Q_1\o$ is a Carter subgroup by 
Lemma \ref{Carter-Cartan-Gen}$(a')$, hence a conjugate of $Q$ by Claim \ref{ClaimTheoRes2}, and 
the maximal nilpotence of $Q$ forces $Q_1\o=Q_1$. Hence Cartan subgroups are definably connected and 
conjugate. 
This completes the proof of Theorem \ref{MainTheoSolvGps}. 
\qed

\begin{Corollary}\label{GConResCartan=Carter}
In a definably connected solvable group definable in an o-minimal structure, 
Cartan subgroups and Carter subgroups coincide.  
\end{Corollary}
\proof
If $Q$ is a Cartan subgroup, then it is definably connected by Theorem \ref{MainTheoSolvGps}, and 
thus a Carter subgroup by Lemma \ref{Carter-Cartan-Gen}$(a')$. 
If $Q$ is a Carter subgroup, then $Q$ is the definably connected component of a Cartan subgroup $\tilde Q$ 
by Lemma \ref{Carter-Cartan-Gen}, and thus $Q=\tilde{Q}$ by Theorem \ref{MainTheoSolvGps}. 
\qed

There are other aspects refining further the structure of definably connected 
solvable groups that we won't follow here, 
but which could be. It includes the already mentioned approach of Cartan/Carter 
subgroups as {\em minimal} abnormal subgroups \cite{Carter61/62, MR1765797}, 
as well as covering properties of nilpotent quotients by Cartan/Carter 
subgroups (see also \cite[\S4-5]{MR2436138}), and also the peculiar 
theory of ``generalized centralizers" of \cite[\S5.3]{MR1765797}. We merely mention 
the most basic covering property, but before that a Frattini argument following 
Theorem \ref{MainTheoSolvGps}. 

\begin{Corollary}\label{CorFrattiniThmRes}
Let $G$ be a group definable in an o-minimal structure, $N$ a definably connected definable normal 
solvable subgroup, and $Q$ a Cartan/Carter subgroup of $N$. Then $G=N_G(Q)N$. 
\end{Corollary}
\proof
By a standard Frattini argument, following the conjugacy in 
Theorem \ref{MainTheoSolvGps}. 
\qed

\begin{Lemma}\label{CarterCoverNilpQuotients}
Let $G$ be a definably connected solvable group definable in an o-minimal structure, $N$ a definable 
normal subgroup such that $G/N$ is nilpotent, and $Q$ a Cartan/Carter subgroup of $G$. Then 
$G=QN$. 
\end{Lemma}
\proof
Suppose $QN<G$. Then $QN/N$ is a definable subgroup of infinite index in the definably connected 
nilpotent group $G/N$. By Lemma \ref{NormalizerCondition}, and since $N_G(QN)$ is the preimage in $G$ 
of $N_{G/N}(QN/N)$, we have thus $QN$ of infinite index in $K:=N_G(QN)$. But $Q$ is a Cartan/Carter 
subgroup of the definably connected solvable group $[QN]\o$, normal in $K$, and thus 
$K=N_K(Q)[QN]\o=N_K(Q)N\o$ by Corollary \ref{CorFrattiniThmRes}. Since $Q$ is a Carter subgroup, 
we get that $QN$ must have finite index in $K$, a contradiction. 
\qed

We note that Lemma \ref{CarterCoverNilpQuotients} always applies with 
$N=F\o(G)$, in view of Fact \ref{FactGenGDefConRes}$(b)$, giving thus 
in particular $G=QF\o(G)$ for any 
definably connected solvable group $G$ and any Cartan/Carter subgroup $Q$ of $G$. 

\section{On Lie groups}\label{SectionLieGroups}

In this section we collect properties needed in the sequel concerning Cartan subgroups 
(in the sense of Chevalley as usual) of Lie groups. These are facts more or less known, 
but because of the different notions of a Cartan subgroup used in the literature 
we will be careful with references. 

By a {\em Lie algebra} we mean a finite dimensional real Lie algebra. 
We are going to make use of the following concepts about Lie algebras: 
subalgebras, commutative, nilpotent, and semisimple Lie algebras 
\cite[I.1.1, I.1.3, I.4.1 and I.6.1]{MR1728312}. 
If $\mathfrak  g$ is a  Lie algebra and $x\in \mathfrak g$, the linear map 
$\ad_x: \mathfrak  g\to  \mathfrak  g: y\mapsto [x,y]$ is called the {\em adjoint  map} of $x$. 
If $ \mathfrak  h$ is a subalgebra of $ \mathfrak  g$, the {\em normalizer} of 
$\mathfrak  h$ in $ \mathfrak  g$ is 
${\mathfrak  n}_{\mathfrak  g}(\mathfrak  h): = \{x\in \mathfrak g \colon \ad_x(\mathfrak  h)\subseteq \mathfrak  h\}$ 
and the {\em centralizer} of $ \mathfrak  h$ in $\mathfrak  g$ is 
${\mathfrak  z}_{\mathfrak  g}(\mathfrak  h):= \{x\in \mathfrak  g\colon [\ad_x]_{|\mathfrak  h}= id_{| \mathfrak  h}\}$. 

\begin{Definition}\label{CSA} 
Let $\mathfrak  g$ be a Lie algebra and $\mathfrak  h$ a subalgebra of $\mathfrak  g$. 
We say that $\mathfrak  h$ is a {\em Cartan subalgebra} of $\mathfrak  g$ if $\mathfrak h$ is 
nilpotent and selfnormalizing in $\mathfrak g$.
\end{Definition}

The two following facts can be found in 
\cite[Theorem 4.1.2]{MR0376938} and \cite[Theorem 4.1.5]{MR0376938} respectively. 

\begin{Fact}\label{Lie0}
Every Lie algebra has  a Cartan subalgebra.
\end{Fact}

\begin{Fact}\label{Lie1} 
Let $ \mathfrak  g$ be a semisimple Lie algebra and $\mathfrak h$ a subalgebra of 
$\mathfrak  g$. Then $\mathfrak  h$ is a Cartan subalgebra of $\mathfrak  g$ if and only if

$(a)$ $\mathfrak h$ is a maximal abelian subalgebra of $\mathfrak  g$, and 

$(b)$ For every $x\in \mathfrak  h$, $\ad_x$ is a semisimple endomorphism of $\mathfrak  g$, i.e., 
$\ad_x$ is diagonalizable over $\C$.
\end{Fact}

By a {\em Lie group} we mean a finite dimensional real Lie group $G$. 
The {\em connected component of the identity} is denoted by $G\o$.
The Lie algebra of $G$  is denoted by  $\mathfrak{L}(G)$. A  connected Lie group $G$
is called  a {\em semisimple Lie group} if $\mathfrak{L}(G)$ is a semisimple
Lie algebra (equivalently, every normal commutative connected immerse subgroup
of $G$ is trivial \cite[Proposition III.9.8.26]{MR1728312}). If $g$ is an
element of a
Lie group $G$, then $\Ad(g): \mathfrak  L (G)\to \mathfrak  L (G)$ denotes the
differential at the identity of $G$ of the map from $G$ to $G$ mapping $h$ to
$ghg^{-1}$, for each $h\in G$. If  $\mathfrak  g$ is the Lie algebra of $G$ and
 $ \mathfrak  h$  a subalgebra of $ \mathfrak  g$, the {\em centralizer} of  
 $\mathfrak h$ in $G$ is $Z_G( \mathfrak  h):=\{g\in G: \Ad(g)(x)=x$ for every
$x\in  \mathfrak  h\}$. 

\begin{Fact}\label{Lie2}
Let G be a connected semisimple Lie group with Lie algebra $\mathfrak  g$, and let $H$ be a subgroup of $G$. 
Then $H$ is a Cartan subgroup of $G$ if and only if $H=Z_G(\mathfrak  h)$ for some Cartan subalgebra 
$\mathfrak  h$ of $\mathfrak  g$. 
Moreover, in this case, $\mathfrak  h$ is $\mathfrak L(H)$.
\end{Fact}
\proof
As $G$ is connected, \cite[Theorem A.4]{MR1367853} implies that $H$ is a Cartan subgroup of $G$ if and only if

(C0)  $H$ is a closed subgroup of $G$;

(C1) $\mathfrak  h(= \mathfrak L (H)$) is a Cartan subalgebra of $\mathfrak  g$, and
 
(C2) $H=C(\mathfrak  h)$.

\noindent
Here $C(\mathfrak  h)$ is defined by a  centralizer-like condition. 
To avoid introducing more notation, instead of properly defining $C(\mathfrak  h)$,  we make use of 
\cite[Lemma I.5]{MR1367853}, which states that $C(\mathfrak  h)=Z_G(\mathfrak  h)$ provided  
$\mathfrak  h$ is reductive in $\mathfrak  g$, which is our case. 
Indeed, $G$ is a semisimple Lie group, so $\mathfrak  g$ is a semisimple Lie algebra, 
hence $\mathfrak  g$ is reductive \cite[Proposition I.6.4.5]{MR1728312}, and
then by 
\cite[Lemma I.4]{MR1367853} every Cartan subalgebra of $\mathfrak  g$ is reductive in $\mathfrak  g$; 
in particular  $\mathfrak  h$ is reductive in  $\mathfrak g$.

For the  converse, we observe that if $H=Z_G(\mathfrak  h)$ for some Cartan
subalgebra $\mathfrak  h$ of $\mathfrak  g$, then $H$ is closed in $G$ and
$\mathfrak L (H)=\mathfrak  h$. Indeed, $H$ is closed by definition of
centralizers, and by \cite[Proposition III.9.3.7]{MR1728312}, 
$\mathfrak L (H)= \mathfrak z_{\mathfrak g}(\mathfrak  h)$. 
Now $\mathfrak  h$ is abelian by Fact \ref{Lie1}, and hence 
$\mathfrak  h\subseteq\mathfrak z_{\mathfrak g}(\mathfrak  h)$. Moreover, if 
$x\in  \mathfrak z_{\mathfrak g}(\mathfrak  h)$, the subalgebra of $\mathfrak g$ generated by $x$ and
$\mathfrak h$ is abelian, so  it must coincide with $\mathfrak h$ by maximality
of $\mathfrak h$, and $x\in \mathfrak h$; hence $\mathfrak h=\mathfrak z_{\mathfrak g}(\mathfrak  h)$. 
We then conclude as above, first  applying Lemma I.5 and then Theorem A.4 from \cite{MR1367853}.
\qed 

\begin{Fact}\label{Lie3}
Let $G$ be a connected semisimple centreless Lie group and $H$ a subgroup of G. 
If $H$ is a Cartan subgroup of $G$, then $H$ is abelian.
\end{Fact}
\proof 
By Fact \ref{Lie2}, $H=Z_G(\mathfrak  h)$ with $\mathfrak  h= \mathfrak L (H)$ 
a Cartan subalgebra of $\mathfrak  g$. By  \cite[Lemma 8, p. 556]{MR0180629} we 
have that $H/Z(G)$ is abelian (see also \cite[Theorem 1.4.1.5]{MR0498999}, 
noting that since $G$ is semisimple 
the general assumption (1.1.5) holds). Hence $H$ is abelian. 
\qed

We note that the assumption $Z(G)=1$ is essential to get the Cartan
subgroup abelian in Fact \ref{Lie3}.  For example $\SL_3(\R)$ has a 
simply-connected double covering with non-abelian Cartan subgroups 
\cite[p.141]{MR568880}, an example which can also occur in the context 
of our Theorem \ref{transfer2} below. 

\begin{Fact}\label{Lie4}
Let $G$ be a connected semisimple Lie group. Then: 
\begin{itemize}
\item[$(a)$]  
There are only finitely many conjugacy classes of Cartan subgroups of $G$. All 
Cartan subgroups of $G$ have the same dimension. 
\item[$(b)$] 
If $H_1$ and $H_2$ are two Cartan subgroups of $G$ with $H_1\o=H_2\o$, 
then $H_1=H_2$. In particular, if $H_1\o$ and $H_2\o$ are conjugate, then $H_1$ and $H_2$ are conjugate 
as well. 
\item[$(c)$]  For any Cartan subgroup $H$ of $G$, the set of elements of 
$H$ contained in a unique conjugate of $H$ is dense in $H$. 
\end{itemize} 
\end{Fact}
\proof
$(a)$. Let $\mathfrak g= \mathfrak L (G)$. Then $\mathfrak g$ is semisimple and it has finitely many 
Cartan subalgebras, say $\mathfrak h_1,\dots,\mathfrak h_s$, such that any Cartan subalgebra 
$\mathfrak h$ of $\mathfrak g$ is conjugate to one of them by an element of $\Ad(G)$, 
i.e., $\Ad(g)(\mathfrak h)=\mathfrak h_i$ for some $i \in \{1,\dots,s\}$ and some $g\in G$ 
(see \cite[Corollary to Lemma 2]{MR0080875} or \cite[Corollary 1.3.1.11]{MR0498999}). 

Next, note that  for every $g$ in $G$ and every (Cartan) subalgebra $\mathfrak h$ of $\mathfrak g$, 
we have $Z_G(\Ad(g)(\mathfrak h))=gZ_G(\mathfrak h)g^{-1}$. For, $h\in Z_G(\Ad(g)(\mathfrak h))$ if and only if 
$\Ad(h)\Ad(g)x=\Ad(g)x$ for every $x\in \mathfrak h$, and the latter is
equivalent to $g^{-1}hg\in Z_G(\mathfrak h)$. Therefore, conjugate Cartan
subalgebras correspond to conjugate centralizers, and by Fact \ref{Lie2} to
conjugate Cartan subgroups.

We prove the second part. By Fact \ref{Lie2} the Lie algebra of a Cartan
subgroup is a Cartan
subalgebra. By \cite[Corollary 4.1.4]{MR0376938} all Cartan
subalgebras have
the same dimension.

$(b)$. It is clear since  $\mathfrak L (H_i)=\mathfrak L (H_i\o)$, for $i=1,2$, and
$H_i=Z_G(\mathfrak L (H_i))$. (Actually, to prove $(b)$ we do not need $G$  to be
semisimple: just consider the $C(\mathfrak L (H_i))$'s of the proof of Fact
\ref{Lie2}, instead of the centralizers.)

$(c)$. We essentially refer to \cite{MR1209133}. Recall, by 
Fact \ref{Lie2} and its proof, that in the 
semisimple case our notion of a Cartan subgroup equals the one used in that paper and 
$C(\mathfrak{h})=Z_G(\mathfrak h)$ for any Cartan subalgebra 
$\mathfrak{h}$ of  $\mathfrak g:= \mathfrak L (G)$. Let $\text{Reg}(G)$ be the set of 
regular elements of $G$, as defined after 
Lemma 1.3 in \cite{MR1209133}. We first show that each element 
$g$ of $\text{Reg}(G)$ lies in a unique Cartan subgroup of 
$G$. Fix $g\in \text{Reg}(G)$. 
By the proof of
\cite[Prop. 1.5]{MR1209133} we have that
$\mathfrak{g}^1(\Ad(g)):=\{x\in \mathfrak{g}:(\exists n\in \mathbb{N})(\Ad(g)-1)^nx=0\}$
is a Cartan subalgebra of $\mathfrak g$ and $g$ belongs to the Cartan
subgroup $Z_G(\mathfrak{g}^1(\Ad(g)))$. To show the uniqueness, let $H$ be a
Cartan subgroup of $G$ containing $g$. By Fact \ref{Lie2},
$H=Z_G(\mathfrak h)$ with $\mathfrak  h= \mathfrak L (H)$
a Cartan subalgebra of $\mathfrak  g$. Since
$g\in Z_G(\mathfrak h)$ we have that
$\mathfrak h\subseteq \mathfrak{g}^1(\Ad(g))$ and hence $\mathfrak{h}=\mathfrak{g}^1(\Ad(g))$ by maximality
of  Cartan subalgebras.
Therefore $H=Z_G(\mathfrak{g}^1(\Ad(g)))$.

Finally, by \cite[Proposition 1.6]{MR1209133}, the subset $\text{Reg}(G)\cap H$ 
is dense in $H$ for all Cartan subgroup $H$ of $G$. 
\qed

For the following, we refer directly to 
\cite[Proposition 5]{MR1871245} and (the proof of) \cite[Lemma 11]{MR1871245} 
respectively. 

\begin{Fact}\label{Lie5}
Let $G$ be a connected Lie group. Then: 
\begin{itemize}
\item[$(a)$]  The union of all Cartan subgroups of $G$ is dense in $G$. 
\item[$(b)$]   For any Cartan subgroup $H$ of $G$, $[H\o]^G$ contains an open subset.
\end{itemize}
\end{Fact}

We finish this section with a remark which, as far as we know, does not seem to have been made before. 
We will show later that all Cartan and Carter subgroups of a group definable in an o-minimal structure are, 
as indicated by Fact \ref{Lie5}$(b)$, weakly generous in the sense of Definition \ref{DefGenerous}$(a)$. 
Our remark is essentially that the stronger notion of generosity of Definition \ref{DefGenerous}$(b)$ 
may be satisfied or not, depending of the Carter subgroups considered, 
and this phenomenon occurs even inside $\SL_2(\R)$. 
Recall that the Cartan subgroups of $\SL_2(\R)$ are, up to conjugacy, the subgroup 
$Q_1$ of diagonal matrices and $Q_2=\SO_2(\R)$. Considering the characteristic polynomial, the two following 
equalities are easily checked: 
$$Q_1^{{\rm SL}_2(\R)}=\{A\in \SL_2(\R): |{\rm tr}(A)|>2\}\cup \{I,-I\}$$ 
$$Q_2^{{\rm SL}_2(\R)}=\{A\in \SL_2(\R): |{\rm tr}(A)|<2\}\cup \{I,-I\}$$ 

\begin{Remark}\label{RemarkSL2Q1GenQ2NotGen}
Let $G=\SL_2(\R)$. Then, according to Definition \ref{DefGenerous}$(b)$:

$(a)$ The Cartan subgroup $Q_1$ of diagonal matrices is generous in $G$. 

$(b)$ The Cartan subgroup $Q_2=\SO_2(\R)$ is not generous in $G$.
\end{Remark}
\proof
$(a)$. 
Fix $a,b\in (0,\frac{1}{13})$ and consider the matrices $A_1=I$, 
$$
A_2=\left( \begin{array}{cc}
            0 & 1 \\
            -1 & 0 \\
           \end{array}
    \right),~
A_3=\left(
\begin{array}{cc}
            a^{-1} & 0 \\
            0 &  a \\
           \end{array}
    \right),\mbox{~and~}
A_4=\left( \begin{array}{cc}
            0 & -b^{-1} \\
            b & 0 \\
           \end{array}
    \right).$$ 
We show that $G=\cup_{i=1}^4 A_iQ_1^G$. Suppose there exists 
$$M=\left( \begin{array}{lr}
            x & y \\
            u & v \\
           \end{array}
    \right)\in G$$ 
with $M\notin \cup_{i=1}^4 A_i Q_1^G$. Since $M\notin A_1Q_1^G\cup
A_2Q_1^G$, we have
$x=\epsilon-v$ and $y=u+\delta$ for some $\epsilon,\delta\in [-2,2]$. 
Since $M \notin A_3Q_1^G$ we have that $|ax+a^{-1}v|=|a(\epsilon-v)+a^{-1}v|\leq 2$, 
so that $v\in [\frac{-2a-a^2\epsilon}{1-a^2},\frac{2a-a^2\epsilon}{1-a^2}]$. 
Since $\epsilon\in [-2,2]$, we deduce that $v\in [\frac{-2a}{1-a},\frac{2a}{1-a}]$. 
Similarly, it follows from $M \notin A_4Q_1^G$ that
$u\in [\frac{-2b}{1-b},\frac{2b}{1-b}]$.
 
Finally, since $a,b<\frac{1}{13}$ we have that $|v|,|u|<\frac{1}{6}$ and $|x|,|y|<2+\frac{1}{6}<3$. 
In particular, $\text{det}(M)=|xv-uy|\leq|x||v|+|u||y|<1$, a contradiction.

$(b)$. 
We show that the family of matrices 
$$M_x=\left( \begin{array}{cc}
            x^2 & x-1 \\
            1 & x^{-1} \\
           \end{array}
    \right)$$ 
with $x>0$ cannot be covered by finitely many translates of $Q_2^G$. It suffices to
prove that for a fixed matrix 
$$A=
\left( \begin{array}{lr}
            a & b \\
            c & d \\
           \end{array}
    \right) \in G$$ 
we have that $\{x\in \R^{>0}: |{\rm tr}(A^{-1}M_x)|>2\}\subseteq \{x\in \R^{>0}:
M_x\notin AQ_2^G\}$ is not bounded. Since ${\rm tr}(A^{-1}M_x)=x^2d-b-c(x-1)+ax^{-1}$ and $x$
is positive, it follows that $|{\rm tr}(A^{-1}M_x)|>2$ if and only if one of the following two conditions holds:
\begin{equation}\label{eq:1}
dx^3-cx^2-(b-c+2)x+a>0
\end{equation}
\begin{equation}\label{eq:2}
dx^3-cx^2-(b-c-2)x+a<0
\end{equation}
It is easy to check that if $d\neq 0$, then either (\ref{eq:1}) or (\ref{eq:2}) is satisfied for large enough $x$. If $d=0$, then $c\neq
0$ (otherwise ${\rm det}(A)=0$) and again the same holds. 
\qed

In Remark \ref{RemarkSL2Q1GenQ2NotGen}, the generous Cartan subgroup is noncompact and 
the nongenerous one is compact. One can then wonder about the various possibilities for generosity 
depending on compactness. 
But considering $Q_1\times Q_2$ in $\SL_2(\R)\times \SL_2(\R)$ one gets from 
Remark \ref{RemarkSL2Q1GenQ2NotGen} a nongenerous and noncompact Cartan subgroup. 
Besides, any compact group is typically covered by a single conjugacy class of compact Cartan subgroups 
by Corollary \ref{CorStructCompactGps}, and these compact Cartan subgroups are in particular generous. 

\section{From Lie groups to definably simple groups}\label{SectionFromLieToDef}

We now return to the context of groups definable in o-minimal structures. 
In the present section we prove the following theorem, essentially transferring 
via Fact \ref{TheoCherlinZilberConjoMin} the results of 
Section \ref{SectionLieGroups} on Lie groups to definably simple groups definable in an 
o-minimal structure. 

\begin{Theorem}\label{transfer}
Let $G$ be a definably simple group definable in an o-minimal structure. 
Then $G$ has definable Cartan subgroups and the following holds.
\begin{itemize}
  \item[$(1)$] $G$ has only finitely many conjugacy classes of Cartan subgroups.
  \item[$(2)$] If $Q_1$ and $Q_2$ are Cartan subgroups of $G$ and  $Q_1\o=Q_2\o$, then  $Q_1=Q_2$.
  \item[$(3)$] Cartan subgroups of $G$ are abelian and have the same dimension.
  \item[$(4)$] If $Q$ is a Cartan subgroup of $G$, then the set of elements of 
  $Q$ contained in a unique conjugate of $Q$ is large in $Q$. In particular, if $a \in Q$, 
  then the set  of elements of $aQ\o$ contained in a unique conjugate of $aQ\o$ is large in $aQ\o$, and 
  $aQ\o$ is weakly generous in $G$.   
  \item[$(5)$] The union of all Cartan subgroups of $G$, which is definable by {\em(1)}, is large in $G$. 
\end{itemize}
\end{Theorem}

Before passing to the proof of Theorem \ref{transfer}, we explain the ``In particular" part of 
item $(4)$. So let $Q$ be a Cartan subgroup such that the set $Q_0$ of elements of $Q$ contained in a 
unique conjugate of $Q$ is large in $Q$. Let $[aQ\o]_0$ be the set of elements of $aQ\o$ 
contained in a unique conjugate of $aQ\o$, for some $a\in Q$. We see easily that 
${Q_0\cap aQ\o}\subseteq [aQ\o]_0$, and since $Q_0$ is large in $Q$ we get that 
${Q_0\cap aQ\o}$ is large in $aQ\o$, as well as $[aQ\o]_0$. Now, since $Q\o \leq N(aQ\o)\leq N(Q\o)$ and 
$\dim(Q\o)=\dim(N(Q\o))$, we get that $\dim([aQ\o]_0)=\dim(N(aQ\o))$, and 
Corollary \ref{CorHGenr=0WH} gives the weak generosity of $aQ\o$. 

We now embark on the proof of Theorem \ref{transfer}, bearing in mind that for item $(4)$ 
we only need to prove the first statement. We first begin with some lemmas. 
By a {\em system of representatives} we mean 
a system of representatives of conjugacy classes of a set of subgroups of a given group. 

\begin{Lemma}\label{lemparam}
Let $\mathcal{M}$ be an o-minimal expansion of an ordered group, $A\subseteq M$ 
a set of parameters containing an element different from $0$, and $G$ a group definable in 
$\mathcal{M}$ over $A$. 
Assume $G$ has, for some $s\in \N$, at least
$s$ non-conjugate Carter subgroups. Then $G$ has at least $s$ non-conjugate
Carter subgroups definable over $A$. In particular, if $G$ has a
finite number of Carter subgroups up to
conjugacy, then in each conjugacy class there exists a Carter subgroup
definable over $A$. 
\end{Lemma}
\proof 
The second part follows easily from the first one. Let $Q_1,\cdots, Q_s$ be
non-conjugate  Carter subgroups of $G$. We denote them by
$Q^{\bar{b}}_{1},\ldots,Q^{\bar{b}}_s$ to stress the fact that they 
are defined over the tuple $\bar{b}$. For each $i=1,\cdots ,s$, let $s_{i}=[N(Q_i):Q_i]$ 
and $r_{i}$ be the nilpotency class of $Q_i$. Consider the set $\Xi$ of tuples $\bar{c}$ satisfying 
the following conditions for each $i$. 

$(1)$ $Q^{\bar{c}}_{i}$ is a nilpotent subgroup of $G$, of nilpotency class $r_{i}$. 

$(2)$ $[N(Q^{\bar{c}}_{i}):Q^{\bar{c}}_{i}]=s_{i}$. 

$(3)$ For any $j=1,\cdots,s$ with $j\neq i$, $Q^{\bar{c}}_{i}$ and $Q^{\bar{c}}_{j}$ are not conjugate. 

$(4)$ $Q^{\bar{c}}_{i}$ is definably connected.

The three first properties are clearly first-order definable. 
The fact that the fourth is also definable is well-known, and for completeness we sketch the proof 
(following Y. Peterzil). 
Let $X\subseteq M^{n+m}$ be a definable set and for each $d\in M^n$ denote by 
$X_d$ the fiber of $X$ over $d$. We have to show that the set 
$\{d\in M^n:X_d \text{ is definably connected}\}$ is definable (here definable 
connectedness is in the topological sense, but by \cite[Proposition 2.12]{MR961362} for a definable 
group the topological notion of definable connectedness coincides with the one generally in use here). 
By the cell decomposition \cite[Thm. III.2.11]{MR1633348}, $X$ is the 
union of definably connected definable sets $C_1,\cdots , C_k$ with the property 
that for each $d\in M^n$ the fiber $(C_i)_d$ is also definably connected. 
Finally, it suffices to note that for each $d\in M^n$ the set 
$X_d=\bigcup_{i=1}^k (C_i)_d$ is definably connected if and only if 
there is an ordering $(C_{i_1})_d,\ldots ,(C_{i_k})_d$ such that 
$\overline{((C_{i_{1}})_d\cup \cdots \cup (C_{i_{j}})_d)}\cap (C_{i_{j+1}})_d\neq \emptyset$ 
or 
$((C_{i_{1}})_d\cup \cdots \cup (C_{i_{j}})_d)\cap\overline{(C_{i_{j+1}})_d}\neq \emptyset$. 

Now the set $\Xi$ is definable, over $A$ since $G$ is, and it is non-empty since it contains $\bar{b}$. 
Since $\mathcal{M}$ expands a group and $A$ contains an element different from $0$, 
the definable closure in $\mathcal{M}$ of $A$ is an elementary substructure of $\mathcal{M}$: 
the theory of $\mathcal{M}$ expanded with a symbol for each element in $A$ has definable
Skolem functions \cite[Chap. 6 \S1(1.1-3)]{MR1633348}, and we may apply the 
Tarski-Vaught test (see also \cite[\S2.3]{MR1773704}). Hence 
there exists a tuple $\bar{c}\in \Xi$ with each coordinate in the definable
closure of $A$. Now $Q^{\bar{c}}_1,\cdots, Q^{\bar{c}}_s$ are
non-conjugate Carter subgroups of $G$, and each can be defined with parameters
in $A$. 
\qed

\begin{Corollary}\label{corparam}
Let $\mathcal{M}$, $A$, and $G$ be as in Lemma {\em \ref{lemparam}}. 
Assume $G$ has a finite number of Cartan subgroups
up to conjugacy. Then, in each conjugacy class there exists a Cartan subgroup
definable over $A$.
\end{Corollary}
\proof 
By Lemma \ref{Carter-Cartan-Gen}, a finite number of conjugacy classes of Cartan
subgroups implies a finite number of conjugacy classes of Carter
subgroups. Hence, by Lemma \ref{lemparam}, there
exists a finite system of representatives of Carter subgroups $Q\o_1,\cdots,Q\o_s$, 
each defined over $A$. Now given any Cartan subgroup $Q$, we have up to conjugacy 
$Q\o=Q\o_i$ for some $i$ by Lemma \ref{Carter-Cartan-Gen}$(a')$, 
and in particular $Q\leq N(Q\o_i)$. Since both
$N(Q\o_i)$ and the finite group $N(Q\o_i)/Q\o_i$ are
definable over $A$, we deduce that $Q$ is definable over $A$ up to conjugacy, as
desired. 
\qed

We will also make use of the following elementary remark, actually valid in any context where 
Lemma \ref{Carter-Cartan-Gen} hold. 

\begin{Remark}\label{RemQioCardinality}
Let $G$ be a group definable in an o-minimal structure such that for 
every pair of Cartan subgroups $Q_1$ and $Q_2$, $Q_1=Q_2$ if and only if $Q\o_1=Q\o_2$. 
Then the cardinality of a system of representatives of Cartan subgroups of $G$ 
equals the cardinality of a system of representatives of Carter subgroups of $G$. 
\end{Remark}

From now on we will use a standard notation from model theory, namely, if $\mathcal{N}_1$ is a 
substructure of $\mathcal{N}_2$ and $X$ is definable in $\mathcal{N}_1$ 
(respectively in $\mathcal{N}_2$ with parameters in $N_1$), then $X(N_2)$ (resp. $X(N_1)$) 
denotes the realization of $X$ in  $\mathcal{N}_2$ (resp. in  $\mathcal{N}_1$). 

\begin{Corollary}\label{cortransfer}
Let $\mathcal{M}$, $A$, and $G$ be as in Lemma {\em \ref{lemparam}}. 
Assume $G$ satisfies properties {\em (1-5)} of Theorem \ref{transfer}. Then:
\begin{itemize}
\item[$(a)$]  
If $\mathcal{N}$ is an elementary substructure of
$\mathcal{M}$ with $A\subseteq N$, then
$G(N)$ also satisfies properties {\em (1-5)}.
\item[$(b)$] 
If $\mathcal{N}$ is an elementary extension of
$\mathcal{M}$, then $G(N)$ also satisfies properties {\em (1-5)}.
\end{itemize}
\end{Corollary}
\proof 
$(a)$. Since $G$ satisfies property (1), it follows from Corollary \ref{corparam} that 
there is a finite system of representatives $Q_1,\cdots, Q_s$ of Cartan subgroups of $G$ defined over $A$. 
Moreover, by Lemma \ref{Carter-Cartan-Gen} and property (2) of $G$ it follows as in Remark \ref{RemQioCardinality} 
that $Q\o_1,\cdots, Q\o_s$ form a system of representatives of  Carter subgroups of $G$ (all defined over $A$). 

We claim  that
\begin{enumerate}
 \item[$(\dagger)$] $Q\o_1(N),\cdots, Q\o_s(N)$  form a system of representatives of Carter subgroups of $G(N)$, and
\item[$(\ddagger)$] $Q_1(N),\cdots, Q_s(N)$  form a system of representatives of Cartan subgroups of $G(N)$.
\end{enumerate}
The claim $(\dagger)$ follows from the definition of a Carter subgroup. Indeed, for each $i\in \{1,\cdots,s\}$, 
since $Q\o_i$ is definably connected, nilpotent, and almost selfnormalizing, $Q\o_i(N)$ satisfies the same properties, 
and is a Carter subgroup of $G(N)$. If $Q\o$ is a Carter subgroup of $G(N)$, then as before $Q\o(M)$ is a
Carter subgroup of $G$, and is $Q\o_i$ for some $i$ up to conjugacy in $G$. Since $\mathcal{N}\preceq \mathcal{M}$, 
$Q\o=Q\o_i(N)$ up to conjugacy in $G(N)$. Similarly, the groups $Q\o_i(N)$ cannot be conjugate 
because the groups $Q\o_i$ are not, proving $(\dagger)$. 

We now show $(\ddagger)$. We first observe: if $R$ is a nilpotent 
definable subgroup of $G(N)$ with $[R\o]^g=Q\o_i(N)$ for some $g\in G(N)$ and 
$i\in \{1,\cdots,s\}$, then $R^g\leq Q_i(N)$. Indeed, 
$Q_i\o=[R\o]^g(M)(=[R\o(M)]^g=[R(M)\o]^g)$. Since 
$[R(M)]^g(=[R^g(M)])$ is nilpotent and $[R(M)^g]\o=[R(M)\o]^g=Q_i\o$ is a Carter subgroup, 
by Lemma \ref{Carter-Cartan-Gen}$(b)$ $R(M)^g$ must be contained in a Cartan 
subgroup which must be $Q_i$ by property (2) of $G$. Therefore $R^g\leq Q_i(N)$, as required. 
Now we deduce $(\ddagger)$ as follows. 
Each $Q_i(N)$ is a Cartan subgroup: by Lemma \ref{Carter-Cartan-Gen} 
there is a Cartan subgroup $Q$ with $Q\o=Q_i(N)\o$ and by the observation above 
we have $Q\leq Q_i(N)$, and $Q=Q_i(N)$ by maximal nilpotence of $Q$. 
It just remains to see that $Q_1(N),\cdots, Q_s(N)$ form a
system of representatives. Let $Q$ be a Cartan subgroup of $G(N)$. By 
Lemma \ref{Carter-Cartan-Gen}$(a')$ $Q\o$ is a Carter subgroup and
then by $(\dagger)$ there exist $g\in G(N)$ and $k\in \{1,\cdots,s\}$
such that $[Q\o]^g=Q\o_k(N)$. Hence $Q^g\leq Q_k(N)$ because of the observation
above, and $Q^g=Q_k(N)$ by maximal nilpotence of $Q$. Finally, observe that 
$Q_1(N),\cdots,Q_s(N)$ cannot be conjugate
in $G(N)$, since $Q_1,\cdots,Q_s$ are not in $G$, proving $(\ddagger)$. 

We now deduce properties (1-5) for $G(N)$ from $(\dagger)$ and $(\ddagger)$. 
Property (1) is exactly $(\ddagger)$. 
For (2), let $R_1$ and $R_2$ be Cartan subgroups of $G(N)$
such that $R_1\o=R_2\o$. By $(\dagger)$, ${R_1\o}^g={R_2\o}^g=Q_i\o(N)$ for some $g\in G(N)$ and some $i$, 
and by the observation in $(\ddagger)$ above we get $R_1^g$, $R_2^g\leq Q_i(N)$, and
an equality by maximal 
nilpotence. In particular $R_1^g=R_2^g$, and $R_1=R_2$. 
Since the dimension in o-minimal structures is 
invariant under elementary substructures, and one considers only definable sets, 
properties (3-5) transfer readily from $G$ to $G(N)$. 

$(b)$. Let $Q_1,\ldots,Q_s$ be a system of representatives of Cartan
subgroups of $G$. By Lemma \ref{Carter-Cartan-Gen} 
and property (2) of $G$ it follows, as in Remark \ref{RemQioCardinality}, that
$Q\o_1,\cdots, Q\o_s$ form a system
of representatives of  Carter subgroups. We first prove that
$Q\o_1(N),\cdots,Q\o_s(N)$ is a
system of representatives of Carter subgroups of $G(N)$. As in $(a)$, we see that
$Q\o_1(N),\cdots,Q\o_s(N)$ are (non-conjugate) Carter
subgroups of $G(N)$. To see that they represent all the conjugacy
classes, suppose there is a Carter subgroup $Q\o$ of $G(N)$ which
is non-conjugate with $Q\o_1(N),\cdots,Q\o_s(N)$. By Corollary \ref{lemparam} we
can assume that $Q\o$ is defined over $M$. Since $Q\o(M)$ is clearly a Carter 
subgroup of $G$, $Q\o(M)^g=Q\o_i$ for some $g\in G$ and some $i$.
Therefore $[Q\o]^g=Q_i\o(N)$, a contradiction.

We next prove that $Q_1(N),\cdots,Q_s(N)$ is a system of representatives of Cartan 
subgroups of $G(N)$. As in $(a)$, it suffices to observe: if $R$ is a nilpotent 
definable subgroup of $G(N)$ with $[R\o]^g=Q\o_i(N)$ for some 
$g\in G(N)$ and $i\in\{1,\ldots,s\}$, then $R^g\leq Q_i(N)$. Indeed, since
$[R\o]^g=Q\o_i(N)$ and $R^g\leq N_{G(N)}(Q\o_i(N))$, $R^g$ is defined over $M$. 
Hence $R^g(M)$ is a definable nilpotent subgroup of $G$ such that
$R^g(M)\o=[R^g]\o(M)=Q\o_i$. Then, by Lemma \ref{Carter-Cartan-Gen}$(b)$ and
property (2) of $G$, $R^g(M)\leq Q_i$. 
In particular $R^g\leq Q_i(N)$, as required. 

Now we can transfer properties (1-5) from $G$ to $G(N)$ as in $(a)$. 
\qed

\noindent
{\bf Proof of Theorem \ref{transfer}.}  
Let $\mathcal{M}$ denote the ground o-minimal structure. 
By Fact \ref{TheoCherlinZilberConjoMin}, there is an $\mathcal{M}$-definable real closed field $R$ 
(with no extra structure) such that $G$ is $\mathcal{M}$-definably isomorphic to a
semialgebraically connected semialgebraically simple semialgebraic group, definable in 
$R$ over the real algebraic numbers $\R_{alg}$. 
By Remark \ref{DimM=DimR}, the dimensions of sets definable in $R$, computed in $\mathcal{M}$ or $R$, are the same. 
Since $\mathcal{M}$-definable bijections 
preserve dimensions, all the conclusions of Theorem \ref{transfer} would then be true if we prove them in this 
semialgebraic group definable in $R$. 
Therefore, replacing $\mathcal{M}$ by $R$, we may suppose that $\mathcal{M}$ is a pure real closed 
field, and that $G=G(M)$ is a semialgebraically connected semialgebraically simple group defined 
over $\R_{alg}$. 

By quantifier elimination $\R_{alg} \preceq R$ and by 
Corollary \ref{cortransfer}$(b)$ it suffices to show our statements for $G(\R_{alg})$. 
By quantifier elimination again, $\R_{alg} \preceq \R$, and by Corollary \ref{cortransfer}$(a)$ 
it now suffices to prove our statements for $G(\R)$. 

Now, we observe that $G(\R)$ is a finite dimensional semisimple centerless connected real Lie group. 
By Facts \ref{Lie0} and \ref{Lie2} it has Cartan subgroups, necessarily definable as usual by 
Lemma \ref{Carter-Cartan-Gen}$(b)$. 
It remains just to notice that all items (1-5) are true in the connected real Lie group $G(\R)$ by 
Facts \ref{Lie3}, \ref{Lie4}, and \ref{Lie5}$(a)$. 
For item (4), we recall that it suffices to prove the first claim, as explained just after the statement 
of Theorem \ref{transfer}. It follows from Fact \ref{Lie4}$(c)$, noticing that a definable subset has maximal dimension 
if and only it has interior \cite[Proposition 2.14]{MR961362}, and thus is dense if and only if it is large. 
\qed

We note that the second claim in Theorem \ref{transfer}(4) could also have been shown using 
Fact \ref{Lie5}$(b)$. 

\section{The semisimple case}\label{SectionSemiSimpleCase}

We now prove a version of Theorem \ref{transfer} for definably connected semisimple groups definable 
in an o-minimal structure. Recall that a definably connected group $G$ is semisimple if $R(G)=Z(G)$ is finite; 
modulo that finite center, $G$ is a direct product of finitely many definably simple groups by 
Fact \ref{TheoCherlinZilberConjoMin}. 

\begin{Theorem}\label{transfer2}
Let $G$ be a definably connected semisimple group definable in an o-minimal structure. 
Then $G$ has definable Cartan subgroups and the following holds.
\begin{itemize}
\item[$(1)$] 
$G$ has only finitely many conjugacy classes of Cartan subgroups.
\item[$(2)$] 
If $Q_1$ and $Q_2$ are Cartan subgroups of $G$ and  $Q_1\o=Q_2\o$, then  $Q_1=Q_2$.
\item[$(3)$] 
If $Q$ is a Cartan subgroup, then $Z(G)\leq Q$, $Q'\leq Z(G)$, and $Q\o\leq Z(Q)$. Furthermore 
all Cartan subgroups have the same dimension.  
\item[$(4)$] 
If $Q$ is a Cartan subgroup of $G$ and $a\in Q$, then the set $[aQ\o]_0$ of elements of 
$aQ\o$ contained in a unique conjugate of $aQ\o$ is large in $aQ\o$, and $aQ\o$ is weakly generous in $G$. 
In addition, if $a_1$ belongs to another Cartan subgroup $Q_1$, then either ${[aQ\o]_0 \cap a_1Q_1\o}= \emptyset$ 
or $aQ\o=a_1Q_1\o$. 
\item[$(5)$] 
The union of all Cartan subgroups of $G$, which is definable by {\em(1)}, is large in $G$. 
In fact, there are finitely many pairwise disjoint definable sets of the form $[aQ\o]_0^G$ with $Q$ a Cartan subgroup 
of $G$ and $a\in Q$, each weakly generic and consisting of pairwise disjoint conjugates of $[aQ\o]_0$, whose union 
is large in $G$. 
\end{itemize}
\end{Theorem}
\proof
Assume first $R(G)=Z(G)=1$. By Fact \ref{TheoCherlinZilberConjoMin}, 
$G=G_1 \times \cdots \times G_n$ where each $G_i$ is an infinite definably simple definable factor. 
Now by Corollary \ref{LemElemCartanDirProd} Cartan subgroups $Q$ of $G$ are exactly of the form 
$$Q=\tilde {Q}_1 \times \cdots \times \tilde {Q}_n$$ 
with $\tilde {Q}_i$ is a Cartan subgroup of $G_i$ for each $i$. In particular $G$ has definable Cartan 
subgroups by Theorem \ref{transfer}. 
Since $Q\o={\tilde {Q}_1\o\times \cdots \times \tilde {Q}_n\o}$ and the dimension is additive, items (1-3) follow easily from 
Theorem \ref{transfer}(1-3). By additivity of the dimension, the first claim in item (4) also transfers 
readily from Theorem \ref{transfer}(4). 
If some element $\alpha$ belongs to $[aQ\o]_0 \cap a_1Q_1\o$, 
for some Cartan subgroups $Q$ and $Q_1$ and some $a\in Q$ and $a_1 \in Q_1$, 
then $Q_1\o\leq C_G\o(\alpha)=Q\o$ by  the commutativity of $Q$ and $Q_1$ and Lemma \ref{FondamentalLemma}, 
and $Q\o=Q_1\o$. 
In particular $aQ\o=\alpha Q\o=\alpha Q_1\o=a_1Q_1\o$, proving item (4). 
For item $(5)$, notice that if some $[aQ\o]_0^G\cap [a_1Q_1\o]_0^G$ is 
non empty in item (4), then $aQ\o=[a_1Q_1\o]^g$ for some $g$ (conjugating in particular $Q\o$ to $Q_1\o$), 
so the finitely many weakly generic definable sets of the form $[aQ\o]_0^G$ are pairwise disjoint and consist of a 
disjoint union of $G$-conjugates of $[aQ\o]_0$. By the largeness of $[aQ\o]_0^G$ in $[aQ\o]^G$ 
provided by Corollary \ref{CorHGenr=0WH} and the largeness of the union of all Cartan subgroups 
provided by Theorem \ref{transfer}(5), the union of all these sets $[aQ\o]_0^G$ is large in $G$, 
proving item $(5)$. 

Assume now just $R(G)=Z(G)$ finite, and let the notation ``$\overline{\phantom{H}}$" denote the quotients by 
$Z(G)$. By the centerless case, all the conclusions of Theorem \ref{transfer2} hold in $\overline{G}$. 
By Lemma \ref{LemElemCartanModZn}, Cartan subgroups of $G$ contain $Z(G)$ and are exactly the preimages 
in $G$ of Cartan subgroups of $G/Z(G)$. In particular, $G$ has definable Cartan subgroups, and we now 
check that they still satisfy (1-5). 

(1) Since $Z(G)$ is contained in each Cartan subgroup, item $(1)$ transfers from the centerless case. 
(2) If $Q_1$ and $Q_2$ are two Cartan subgroups of $G$ with $Q_1\o=Q_2\o$, then 
$\overline{Q_i\o}=[\overline{Q_i}]\o$ and $\overline{Q_1}=\overline{Q_2}$ by $(2)$ in $\overline{G}$, 
giving $Q_1=Q_2$. 
(3) By the centerless case $\overline{Q}$ is abelian, and thus $Q'\leq Z(G)$. 
In particular $[Q,Q\o]$ is in the finite center $Z(G)$, but since $[Q,Q\o]$ is definable and definably 
connected by \cite[Corollary 6.5]{BaroJalOteroCommutators} we get $[Q,Q\o]=1$, proving the first claim of $(3)$. 
Since the natural (and definable) projection from $G$ onto $\overline{G}$ has finite fibers one gets by 
axioms A2-3 of the dimension that $\dim(\overline{Q})=\dim(Q)$, transferring also from $\overline{G}$ to $G$ 
the second claim of (3). 
(4) Let $Q$ and $Q_1$ be two Cartan subgroups, $a\in Q$ and $a_1\in Q_1$. If some element $\alpha$ belongs 
to $[aQ\o]_0 \cap a_1Q_1\o$, one sees as in the centerless case, still using Lemma \ref{FondamentalLemma} 
but now the fact that $Q\o\leq Z(Q)$ and $Q_1\o \leq Z(Q_1)$, that $aQ\o = a_1Q_1\o$. 
We now show that $[aQ\o]_0$ is large in $aQ\o$. For that purpose, first notice that $[aQ\o]_0$ is exactly the set of 
elements of $aQ\o$ contained in finitely many conjugates of $aQ\o$: for, if $\alpha$ is in $aQ\o$ and in only finitely many 
of its conjugates, say ${(aQ\o)}^{g_1}, \cdots {(aQ\o)}^{g_k}$, then as above 
Lemma \ref{FondamentalLemma} yields 
$Q\o=C\o(\alpha)$, and $aQ\o={(aQ\o)}^{g_1}= \cdots ={(aQ\o)}^{g_k}$. 
For the largeness of $[aQ\o]_0$ in $aQ\o$, it suffices as in item (3) to show that 
$[aQ\o]_0$ contains the preimage of the set of elements $\overline{\alpha}$ of 
$\overline{aQ\o}$ contained in a unique $\overline{G}$-conjugate of 
$\overline{aQ\o}$. So assume towards a contradiction that there exists 
an element $\alpha$ in $aQ\o$, in 
infinitely many $G$-conjugates of $aQ\o$ but such that $\overline{\alpha}$ is in a unique conjugate of 
$\overline{aQ\o}$. Now for $g$ varying in infinitely many cosets of $N(aQ\o)$, 
and in particular in infinitely many cosets of $N\o(aQ\o)=N\o(Q\o)=Q\o$, 
we have $aZ(G)Q\o=[aZ(G)Q\o]^g$. 
But such elements $g$ must normalize the subgroup $Z(G)Q\o$, and in particular 
$[Z(G)Q\o]\o=Q\o$, and hence cannot vary in infinitely many cosets of $Q\o$. This contradiction proves that 
$[aQ\o]_0$ is large in $aQ\o$, and the weak generosity of $aQ\o$ in $G$ follows as usual with 
Corollary \ref{CorHGenr=0WH}. 
(5) Using the projection from $G$ to $\overline{G}$, the non weak genericity of the complement of the union 
of all Cartan subgroups passes from $\overline{G}$ to $G$, and thus the union of all Cartan subgroups of $G$ 
is large in $G$. Then all other claims of item $(5)$ follow as in the case $Z(G)=1$. 
\qed

In Theorem \ref{transfer2}(3) Cartan subgroups need not be abelian outside of the centerless case, since 
the simply-connected double covering of $\SL_3(\R)$ with non-abelian Cartan subgroups 
mentioned after Fact \ref{Lie3} is definable in $\R$. 
The following question then arises naturally. 

\begin{Question}
Let $G$ be a definably connected semisimple group definable in an o-minimal structure, 
and $Q$ a Cartan subgroup of $G$. When is it the case that $Q$ is abelian? That $Q=Q\o Z(G)$?
\end{Question}

For Carter subgroups, one gets the following corollary of Theorem \ref{transfer2}. 

\begin{Corollary}
Let $G$ be a definably connected semisimple group definable in an o-minimal structure. 
Then $G$ has finitely many conjugacy classes of Carter subgroups. 
Each Carter subgroup $Q\o$ is abelian and weakly generous in the following strong sense: 
the set of elements of $Q\o$ contained in a unique conjugate of $Q\o$ is large in $Q\o$ 
and weakly generous in $G$. 
\end{Corollary}
\proof
We know by Lemma \ref{Carter-Cartan-Gen} that Carter subgroups are exactly the definably connected 
components $Q\o$ of Cartan subgroups $Q$ of $G$. 
In particular item $(3)$ of Theorem \ref{transfer2} shows that 
$Q\o\leq Z(Q)$, and $Q\o$ is abelian. The other conclusions follow immediately from 
items $(1)$ and $(4)$ in Theorem \ref{transfer2}. 
\qed

Before moving to more general situations, we make a few additional remarks about the semisimple case. 
We first mention a general result on control of fusion, reminiscent from \cite[Corollary 2.12]{DeloroJaligotI} 
in the finite Morley rank case. 

\begin{Lemma}[Control of fusion]\label{LemmaControFusion}
Let $G$ be a group definable in an o-minimal structure, $Q$ a Cartan subgroup of $G$, $X$ and $Y$ 
two $G$-conjugate subsets of $C(Q\o)$ such that $C\o(Y)$ has a single 
conjugacy class of Carter subgroups. Then $Y=X^g$ for some $g$ in $N(Q\o)$. 
\end{Lemma}
\proof
Let $g$ in $G$ be such that $Y=X^g$. Then $C\o(Y)=C\o(X)^g$ contains both $Q\o$ and ${Q\o}^g$, 
so our assumption forces that ${[Q\o]}^{g\gamma}={Q\o}$ for some $\gamma$ in $C\o(Y)$. 
Now $g\gamma$ normalizes $Q\o$ and $X^{g\gamma}=Y^{\gamma}=Y$. 
\qed

\begin{Lemma}\label{Q=FNQ}
Let $G$ be a definably connected semisimple group $G$ definable in an o-minimal structure and 
$Q$ a Cartan subgroup of $G$. Then $Q=F(N_G(Q\o))$. 
\end{Lemma}
\proof
Any definable nilpotent subgroup containing the Carter subgroup $Q\o$ is a finite extension 
of it by Lemma \ref{NormalizerCondition}, and hence is in $N_G(Q\o)$. By Theorem \ref{transfer2}(2), 
there is a unique maximal one. This proves that $Q\trianglelefteq N_G(Q\o)$. Hence $Q\leq F(N_G(Q\o))$, and 
in fact there is equality by maximal nilpotence of $Q$. 
\qed

With Lemma \ref{LemmaControFusion}, we can rephrase the last part of Theorem \ref{transfer2}(4). 

\begin{Corollary}\label{CorControlFusionSemisimple}
Let $G$ be a definably connected semisimple group definable in an o-minimal structure and 
$Q$ a Cartan subgroup of $G$. If $a_1$ and $a_2$ are two $G$-conjugate elements of $Q$ 
such that $a_i\in [a_iQ\o]_0$ as in Theorem \ref{transfer2}{\em (4)} for $i=1$ and $2$, 
then $a_1Q\o$ and $a_2Q\o$ are $N(Q)$-conjugate.  
\end{Corollary}
\proof
By Theorem \ref{transfer2}(3), $a_i\in C(Q\o)$ for each $i$, and by Lemma \ref{FondamentalLemma} 
$Q\o=C\o(a_1)=C\o(a_2)$. Lemma \ref{LemmaControFusion} implies then that $a_2={a_1}^g$ for some 
$g$ in $N(Q\o)$. But since $Q\trianglelefteq N_G(Q\o)$ by Lemma \ref{Q=FNQ}, $g\in N_G(Q)$. 
\qed

As just seen in Corollary \ref{CorControlFusionSemisimple}, if $Q$ is a Cartan 
subgroup of a definably connected semisimple group $G$ definable in an o-minimal structure, then 
$$N_G(Q)=N_G(Q\o).$$ 
Now the finite group 
$W(G,Q):=N_G(Q)/Q=N_G(Q\o)/Q$ 
can naturally be called the {\em Weyl group relative to $Q$}, or, equivalently, {\em relative to $Q\o$}. 
If $G$ is definably simple, then one has the two alternatives at the end of 
Fact \ref{TheoCherlinZilberConjoMin}. In the first case $G$ is essentially a simple algebraic group over an 
algebraically closed field (of characteristic $0$). It is well known in this case that there is only one conjugacy of 
Cartan subgroups, the maximal (algebraic and split) tori which are also Carter subgroups (by divisibility). 
Then there is only one relative Weyl group, and their classification is provided 
by the classification of the simple algebraic groups. In the second alternative at the end of 
Fact \ref{TheoCherlinZilberConjoMin}, the group is essentially a simple real Lie group, and again 
the Weyl groups relative to the various Cartan subgroups, corresponding to the various split or non-split tori, 
are classified in this case. For a general definably connected semisimple ambient group $G$, the structure 
of the Weyl groups is inherited from that of the definably simple factors of $G/R(G)$, as we will see in 
Section \ref{SectionFinalComments}. 

Theorem \ref{transfer2}(5) equips any definably connected semisimple group with some kind of a 
partition into finitely many canonical ``generic types". We finish this section by counting them precisely. 

\begin{Remark}\label{RemarknG}
The number $n(G)$ of weakly generic definable sets of the form $[aQ\o]_0^G$ as in 
Theorem \ref{transfer2}$(5)$ is clearly bounded by the sum $\Sigma_{Q\in \mathcal{Q}} |Q/Q\o|$ 
where $\mathcal Q$ is a system of representatives of the set of Cartan subgroups of $G$. 
But it might happen in Theorem \ref{transfer2}$(4)$ that two distinct sets of the form $aQ\o$ and 
$a'Q\o$, for $a$ and $a'$ in a common Cartan subgroup $Q$, are conjugate by the action of the 
Weyl group $W(G,Q)=N_G(Q)/Q$. If one denotes 
by $\sim_Q$ the equivalence relation on $Q/Q\o$ by the action of $W(G,Q)$ 
naturally induced by conjugation on $Q/Q\o$, then one sees indeed with Corollary \ref{CorControlFusionSemisimple} that  
$$n(G)=\Sigma_{Q\in {\mathcal Q}}|[Q/Q\o]_{\sim_Q}|.$$ 
\end{Remark}

\section{The general case}\label{SectionGeneralCase}

We now analyze the general case of a group definable in an o-minimal structure. 
As far as possible, we will restrict ourselves to definably connected groups only when necessary. 
We start by lifting Carter subgroups. 

\begin{Lemma}\label{LiftCarterGeneral}
Let $G$ be a group definable in an o-minimal structure, and $N$ a definable normal subgroup of $G$ 
such that $N\o$ is solvable. Then Carter subgroups of $G/N$ are exactly of the form 
$QN/N$ for $Q$ a Carter subgroup of $G$. 
\end{Lemma}
\proof
We may use the notation ``$\overline{\phantom{H}}$" to denote the quotients by $N$. 
Let $Q$ be a Carter subgroup of $G$. Then $Q$ is also a Carter subgroup of the definable subgroup $QN$. 
The preimage in $G$ of $N_{\overline{G}}(\overline{Q})$ normalizes $[QN]\o=QN\o$, and thus is contained in 
$N_G(Q)N$ by Corollary \ref{CorFrattiniThmRes}. Hence $\overline{Q}$, which is definable and definably connected, 
must have finite index in its normalizer in $\overline{G}$, and is thus a Carter subgroup of $\overline{G}$. 
Conversely, let $X/N$ be a Carter subgroup of $\overline{G}$ for some subgroup $X$ of $G$ containing $N$. 
Since $X/N$ is definable, $X$ must be definable. By Theorem \ref{MainTheoSolvGps}, 
$X\o$ has a Carter subgroup $Q$, and of course $Q$ must also be a Carter subgroup of $X$. 
Since $X=X\o N$ and $X\o=Q(X\o \cap N)$ by Lemma \ref{CarterCoverNilpQuotients}, we get that 
$X=QN$. Since $QN/N$ is a Carter subgroup $\overline{G}$, we get that $QN$ has finite 
index in $N_G(QN)$. Since $N_G(Q)\leq N_G(QN)$ and $Q$ has finite index in $N_{QN}(Q)$, we get 
that $Q$ has finite index in $N_G(Q)$. Hence $X=QN$ for a Carter subgroup $Q$ of $G$. 
\qed

The following special case of Lemma \ref{LiftCarterGeneral} with $N=R\o(G)$ is of 
major interest, and for the rest of the paper one should bear in mind that 
$$R\o(G)=R\o(G\o).$$ 

\begin{Corollary}\label{CorLiftCarterGeneral}
Let $G$ be a group definable in an o-minimal structure. Then Carter subgroups of $G/R\o(G)$ are exactly 
of the form $QR\o(G)/R\o(G)$ for $Q$ a Carter subgroup of $G$. 
\end{Corollary}

At this stage, we can prove our general Theorem \ref{MainTheo} giving the existence, the definability, and 
the finiteness of the set of conjugacy classes of Cartan subgroups in an arbitrary group definable in 
an o-minimal structure. 

\bigskip
\noindent
{\bf Proof of Theorem \ref{MainTheo}.}  
Let $G$ be our arbitrary group definable in an arbitrary o-minimal structure. The quotient 
$G\o/R\o(G)$ is semisimple by Fact \ref{ActionConOnFinite}, 
and has Carter subgroups by Theorem \ref{transfer2}. 
Hence $G\o$ has Carter subgroups 
by Corollary \ref{CorLiftCarterGeneral}. This takes care 
of the existence of Carter subgroups of $G\o$, and of course of $G$ as well. 
Now $G$ has Cartan subgroups by Lemma \ref{Carter-Cartan-Gen}. 
Their definability is automatic as usual in view of Lemma \ref{Carter-Cartan-Gen}$(a')$. 
To prove that Cartan subgroups fall into only finitely many conjugacy classes, it suffices by 
Lemma \ref{Carter-Cartan-Gen}$(a')$ to prove it for Carter subgroups. We may then assume $G$ 
definably connected. Now groups of the form $QR\o(G)/R\o(G)$, for $Q$ a Carter subgroup of $G$, are 
Carter subgroups of the semisimple quotient $G/R\o(G)$. By Theorem \ref{transfer2}(1), there are only 
finitely many $G/R\o(G)$-conjugacy classes of groups of the form $QR\o(G)/R\o(G)$, and thus only 
finitely many $G$-conjugacy classes of groups of the form $QR\o(G)$. Replacing $G$ by 
such a $QR\o(G)$, we may thus assume $G$ definably connected and solvable. 
But now in $G$ there is only one conjugacy class of Carter subgroups by Theorem \ref{MainTheoSolvGps}. 
This completes our proof of Theorem \ref{MainTheo}. 
\qed

We mention the following form of a Frattini Argument as a consequence of Theorem \ref{MainTheo}. 

\begin{Corollary}
Let $G$ be a definably connected group definable in an o-minimal structure and $N$ a definable normal 
subgroup of $G$. Then $G=N_G\o(Q)N\o$ for any Cartan subgroup $Q$ of $N$. 
\end{Corollary}
\proof
Clearly, for any element $g$ of $G$, $Q^g$ is a Cartan subgroup of $N$. On the other hand, the 
set $\mathcal Q$ of conjugacy classes of Cartan subgroups of $N$ is finite by Theorem \ref{MainTheo}, 
and the action of $G$ on $N$ by conjugation naturally induces a definable action on the finite set $\mathcal Q$. 
Since $G$ is definably connected, Fact \ref{ActionConOnFinite}$(a)$ shows that this action must be trivial. 
Hence, for any $g$ in $G$, $Q^g$ is indeed in the same $N$-conjugacy class as $Q$, i.e., 
$Q^g=Q^h$ for some $h\in N$; in particular $g=gh^{-1}h \in N_G(Q)N$. 
Hence $G=N_G(Q)N$, and in fact $G=N_G\o(Q)N\o$ by definable connectedness. 
\qed

We shall now inspect case by case what survives of Theorem \ref{transfer2}(2-5) in the general case. 
We first consider Theorem \ref{transfer2}(2). 

\begin{Theorem}\label{TheoGen(2)}
Let $G$ be a definably connected group definable in an o-minimal structure and $Q$ a Cartan subgroup 
of $G$. Then there is a  unique (definable) subgroup $K_Q$ of $G$ containing $R\o(G)$ 
such that $K_Q/R\o(G)$ is the unique Cartan subgroup of $G/R\o(G)$ containing ${Q\o R\o(G)}/R\o(G)$. 
Moreover, $QR(G)\leq K_Q$ and 
$$Q=F(N_{K_Q}(Q\o))=C_{K_Q}(Q\o) Q\o=C_G(Q\o) Q\o.$$
\end{Theorem}
\proof
By Corollary \ref{CorLiftCarterGeneral}, the group ${Q\o R\o(G)}/R\o(G)$ is a Carter subgroup of the 
semisimple quotient $G/R\o(G)$. By Theorem \ref{transfer2}(2), it is contained in a unique Cartan subgroup, 
of the form $K/R\o(G)$ for some subgroup $K$ containing $R\o(G)$ and necessarily definable by 
Lemma \ref{Carter-Cartan-Gen}$(a')$. We will show that $K_Q=K$ satisfies all our claims. 
Since ${QR\o(G)}/R\o(G)$ is nilpotent and contains the Carter subgroup 
${Q\o R\o(G)}/R\o(G)$, we have $QR\o(G)\leq K$. Since $R(G)/R\o(G)$ is the center of $G/R\o(G)$, it is 
contained in $K/R\o(G)$ by Lemma \ref{LemElemCartanModZn}$(a)$, and thus $R(G)\leq K$. 
Hence, $QR(G)\leq K$. 

To prove our last equalities, we first show that $F(N_K(Q\o))=C_K(Q\o)Q\o$. 
Since $Q\o=F\o(N_K(Q\o))$ by Lemma \ref{NormalizerCondition}, 
the inclusion from left to right follows from Fact \ref{StructNilpGps}. For the reverse inclusion, 
notice that ${C_K(Q\o)Q\o}$ is normal in $N_K(Q\o)$. Since Cartan subgroups of 
$G/R\o(G)$ are nilpotent in two steps by Theorem \ref{transfer2}(3), the second term of the 
descending central series of ${C_K(Q\o)Q\o}$ is in $R\o(G)$, and thus in 
$Q\o$ because $Q\o$ is selfnormalizing in $Q\o R\o(G)$ by 
Theorem \ref{MainTheoSolvGps}. By keeping taking descending central series and using the nilpotency 
of $Q\o$, we then see that ${C_K(Q\o)Q\o}$ is nilpotent, and thus in 
$F(N_K(Q\o))$ by normality in $N_K(Q\o)$. 

Since $C_G(Q\o)\leq K$, clearly by considering its image modulo $R\o(G)$, our last equality is true. 
Finally, $Q=C_Q(Q\o)Q\o$ by Fact \ref{StructNilpGps}, and thus $Q\leq C_K(Q\o)Q\o=F(N_K(Q\o))$. 
Now the maximal nilpotence of $Q$ forces $Q=F(N_K(Q\o))$, and our proof is complete. 
\qed

With Theorem \ref{TheoGen(2)} one readily gets the analog of Theorem \ref{transfer2}(2). Of course 
definable connectedness is a necessary assumption here, since a finite group may have several 
Cartan subgroups. 

\begin{Corollary}\label{CorQ1oQ2oGenCase}
Let $G$ be a definably connected group definable in an o-minimal structure, $Q_1$ and $Q_2$ 
two Cartan subgroups. If $Q\o_1=Q\o_2$, then $Q_1=Q_2$. 
\end{Corollary}

We also get that $QR\o(G)$ is normal in $K_Q$, and actually has a quite stronger 
uniqueness property in $K_Q$. 

\begin{Corollary}\label{CorTheoGen(2)}
Same assumptions and notation as in Theorem \ref{TheoGen(2)}. Then $[K_Q]\o={Q\o R\o(G)}$ and 
$QR\o(G)$ is invariant under any automorphism of $K_Q$ leaving $[K_Q]\o$ invariant. 
\end{Corollary}
\proof
The first equality comes from Lemma \ref{LiftCarterGeneral}. 

Let $\sigma$ be an arbitrary automorphism of $K_Q$ leaving $[K_Q]\o$ invariant. 
Since $Q\o$ is a Cartan subgroup of $[K_Q]\o$ by Corollary \ref{GConResCartan=Carter}, 
its image by $\sigma$ is also a Cartan subgroup of $[K_Q]\o$, 
and with Theorem \ref{MainTheoSolvGps} one gets $[Q\o]^{\sigma}=[Q\o]^k$ for some $k$ in $[K_Q]\o$. 
Since $QR\o(G)$ is normalized by $k$, we can thus assume that $\sigma$ leaves $Q\o$ invariant. 
But now $\sigma$ leaves $F(N_{K_Q}(Q\o))$ invariant. Hence by 
Theorem \ref{TheoGen(2)} $Q$ is left invariant by $\sigma$, and thus $\sigma$ leaves 
$Q[K_Q]\o=QR\o(G)$ invariant.  
\qed

The main question we are facing with at this stage is the following. 

\begin{Question}\label{QuestLiftCartans}
Is it the case, in Theorem \ref{TheoGen(2)}, that $K_Q=QR\o(G)$?
\end{Question}

Question \ref{QuestLiftCartans} has a priori stronger forms, which are indeed equivalent 
as the following lemma shows. 

\begin{Lemma}\label{LemKQEqualsPRQR}
Under the assumptions and notation of Theorem \ref{TheoGen(2)}, the following are equivalent: 

$(a)$ $K_Q=PR\o(G)$ for some Cartan subgroup $P$ of $G$

$(b)$ $K_Q=PR\o(G)$ for any Cartan subgroup $P$ of $K_Q$. 
\end{Lemma}
\proof
Assume $K_Q=P_1R\o(G)$ for some Cartan subgroup $P_1$ of $G$, and suppose 
$P_2$ is a Cartan subgroup of $K_Q$. Then $P\o_1$ and $P\o_2$ are Carter subgroups of 
$[K_Q]\o$ by Lemma \ref{Carter-Cartan-Gen}$(a')$. Since they are $[K_Q]\o$-conjugate by 
Theorem \ref{MainTheoSolvGps}, we may assume $P\o_1=P\o_2$ up to conjugacy. Now applying 
Theorem \ref{TheoGen(2)} with the Cartan subgroup $P_1$, or just Corollary \ref{CorQ1oQ2oGenCase}, 
we see that $P_1=P_2$ up to conjugacy, and thus $K_Q=P_2R\o(G)$. 

Conversely, suppose $K_Q=PR\o(G)$ for any Cartan subgroup $P$ of $K_Q$. This applies in particular 
to the Cartan subgroup $Q$ of $G$. 
\qed

By the usual Frattini Argument following the conjugacy of Cartan/Carter subgroups in 
$[K_Q]\o$, we have that $K_Q=\hat{Q}R\o(G)$ where 
$$\hat{Q}=N_{K_Q}(Q\o).$$
The subgroup $\hat{Q}$ is solvable and nilpotent-by-finite, and with the selfnormalization property 
of $Q\o$ in the definably connected solvable group $Q\o R\o(G)$ one sees easily that 
$\hat{Q}/Q \simeq {K_Q/{(QR\o(G))}}$. 
Hence Question \ref{QuestLiftCartans} is equivalent to proving that the finite quotient 
$\hat{Q}/Q$ is trivial. 

Retaining all the notation introduced so far, Theorem \ref{transfer2}(3) takes the following form 
for a general definably connected group. 

\begin{Theorem}\label{TheoGen(3)}
Same assumptions and notation as in Theorem \ref{TheoGen(2)}. Then 
$[K_Q]'\leq R(G)$, and $[\hat{Q},[\hat{Q}]']\leq {Q\o\cap R\o(G)}$ where 
$\hat{Q}=N_{K_Q}(Q\o)$. 
\end{Theorem}
\proof
By Theorem \ref{transfer2}(3), $[K_Q]'\leq R(G)$ and $[K_Q,[K_Q]']\leq R\o(G)$. 
The second inclusion shows in particular that 
$[\hat{Q},[\hat{Q}]']\leq R\o(G)$, and since $Q\o$ is selfnormalizing in 
$Q\o R\o(G)$ by Theorem \ref{MainTheoSolvGps}, we get inclusion in $Q\o$ as well. 
\qed

We now consider Theorem \ref{transfer2}(4) and give its most general form in the general case 
(working in particular without any assumption of definable connectedness of the ambient group). 

\begin{Theorem}\label{TheoGen(4)}
Let $G$ be a group definable in an o-minimal structure, $Q$ a Cartan subgroup of $G$ and $a\in Q$. 
Then $aQ\o$ is weakly generous in $G$. In fact, the set of elements of $aQ\o$ contained in a unique conjugate 
of $aQ\o$ is large in $aQ\o$. Furthermore, if $G$ is definably connected, then 
the set of elements of $Q$ contained in a unique conjugate of $Q$ is large in $Q$.  
\end{Theorem}
\proof
We first prove that the set of elements of $Q\o$ contained in a unique $G$-conjugate of $Q\o$ 
is large in $Q\o$. For that purpose, it suffices by Corollary \ref{CorDescriptionQwgenCarter} to 
show that the set of elements of $Q\o$ contained in only finitely many $G$-conjugates of $Q\o$ is large 
in $Q\o$. Assume towards a contradiction that the set $Q_\infty$ of elements of $Q\o$ contained in infinitely 
many $G$-conjugates of $Q\o$ is weakly generic in $Q\o$. By Theorem \ref{MainTheoSolvGps}, we may restrict 
$Q_\infty$ to the subset of elements contained in a unique $Q\o R\o(G)$-conjugate of $Q\o$, and still 
have a weakly generic subset of $Q\o$. 
Now $Q_\infty$ must have a weakly generic image in $Q\o$ modulo $R\o(G)$. 
By Theorem \ref{transfer2}(4), we must then find an element 
$x\in Q_\infty$ which, modulo $R\o(G)$, 
is in a unique conjugate of $Q\o$. Then we have infinitely many Carter subgroups of 
$Q\o R\o(G)$ passing through $x$, a contradiction since they are all $Q\o R\o(G)$-conjugate by 
Theorem \ref{MainTheoSolvGps}. 

We now consider the full Cartan subgroup $Q$, and an arbitrary element $a$ in $Q$. 
For the weak generosity of $aQ\o$ in $G$, it suffices to use our general Corollary \ref{CorHGenCosetwHGen}. 
Indeed, by Corollary \ref{CorHGenr=0WH}, 
it suffices to show the stronger property that the set of elements of $aQ\o$ in a unique 
conjugate of $aQ\o$ is large in $aQ\o$. Assume towards a contradiction that the set $X$ of elements of 
$aQ\o$ in at least two distinct conjugates of $aQ\o$ is weakly generic in $aQ\o$. If $n$ is the order 
of $a$ modulo $Q\o$, then the set of $n$-th powers of elements of $X$ would be weakly generic 
in $Q\o$ by Corollary \ref{CorLemPhi(X)Gen}. 
Hence by the preceding paragraph one would find an element $x$ in $X$ such that $x^n$ 
is in a unique conjugate of $Q\o$. This is a contradiction as usual 
since $xQ\o$ must then be the unique conjugate of $aQ\o$ containing $x$. 

We now prove our last claim about $Q$ when $G$ is definably connected. 
Assume towards a contradiction that the set 
$X$ of elements in $Q$ and in at least two distinct conjugates of $Q$ is weakly generic in $Q$. Then it 
should meet one of the cosets $aQ\o$ of $Q\o$ in $Q$ in a weakly generic subset, say $X'$. 
By Corollary \ref{CorLemPhi(X)Gen} again, one finds an element $x$ in $X'$ such that $x^{|Q/Q\o|}$ 
is in a unique conjugate of $Q\o$. Now all the conjugates of $Q$ passing through $x$ should have the same 
definably connected component, and thus are $N_G(Q\o)$-conjugate. Then they are all equal by 
Corollary \ref{CorQ1oQ2oGenCase}, a contradiction. 
\qed

In case Question \ref{QuestLiftCartans} fails, we unfortunately found no way of proving 
Theorem \ref{TheoGen(4)} for $a$ in $\hat{Q} \setminus Q$. 
Besides, our method for proving the weak generosity of $aQ\o$ in $G$ does not seem to be 
appropriate for attacking the following more refined question. 

\begin{Question}\label{QuestaQolqrgeGen}
Let $G$, $Q$, and $a$ be as in Theorem \ref{TheoGen(4)}, with $G$ definably connected and 
such that, modulo $R\o(G$), $a$ is in a unique conjugate of $aQ\o$. 

$(a)$ Is it the case that $[aQ\o]^{R\o(G)}$ is large in $aQ\o R\o(G)$? 

$(b)$ Same question, with $a$ in $\hat{Q}$ instead of $a$ in $Q$?
\end{Question}

By Theorem \ref{TheoGen(4)}, the union of Cartan subgroups of a group definable in an o-minimal structure 
must be weakly generic, but the much stronger statement of Theorem \ref{transfer2}(5) now becomes a 
definite question. 

\begin{Question}\label{QuestQlarges}
Let $G$ be a definably connected group definable in an o-minimal structure. Is it the case that the 
union of its Cartan subgroups forms a large subset? 
\end{Question}

We now prove that Question \ref{QuestQlarges} can be seen on top of both 
Questions \ref{QuestLiftCartans} and \ref{QuestaQolqrgeGen}. 

\begin{Proposition}\label{PropLargeImpliesAll}
Let $G$ be a definably connected group definable in an o-minimal structure 
whose Cartan subgroups form a large subset. Then 
\begin{itemize}
\item[$(a)$] 
Cartan subgroups of $G/R\o(G)$ are exactly of the form $QR\o(G)/R\o(G)$ with $Q$ a Cartan 
subgroup of $G$. 
\item[$(b)$] 
For every Cartan subgroup $Q$ and $a$ in $Q$ such that, modulo $R\o(G)$, $a$ is in a unique 
conjugate of $aQ\o$, $[aQ\o]^{R\o(G)}$ is large in $aQ\o R\o(G)$. 
\end{itemize}
\end{Proposition}
\proof
$(a)$. Assume towards a contradiction that for some Cartan subgroup $Q$, and with the previously used 
notation, we have $QR\o(G)<K_Q$. 
Let $\overline{B}$ be the large subset of ${(K_Q/R\o(G))}\setminus {(QR\o(G)/R\o(G))}$ then provided 
by Theorem \ref{transfer2}(4), and $B$ its pull back in $G$. 
By additivity of the dimension, $B^G$ must be weakly generic in $G$. 
Now the largeness of the set of Cartan subgroups forces the existence of an element $g$ in 
$B \cap P$ for some Cartan subgroup $P$ of $G$. 
Let $\overline{g}$ denote the image of $g$ in $G/R\o(G)$. 
We have $g\in {K_Q\setminus QR\o(G)}$, and 
$C\o(\overline{g})=Q\o R\o(G) / R\o(G)$ 
by considering the structure of Cartan subgroups in the semisimple 
quotient $G/R\o(G)$ and the uniqueness property of $\overline{g}$. 
By Lemma \ref{LiftCarterGeneral} the group $P\o$, modulo $R\o(G)$, is a Carter subgroup of $G/R\o(G)$. 
Now $P$, modulo $R\o(G)$, is included in a Cartan subgroup of $G/R\o(G)$, and its definably connected 
component centralizes $\overline{g}$ by Theorem \ref{transfer2}$(3)$. We then get 
$P\o R\o(G)/R\o(G) \leq C\o(\overline{g}) ={Q\o R\o(G)/R\o(G)}$, and actually equality since the first group 
is a Carter subgroup. Hence $P\o R\o(G)=Q\o R\o(G)$ and Theorem  \ref{TheoGen(2)} yields $P\leq K_Q$. 
Since $Q\o$ and $P\o$ are conjugate in $Q\o R\o(G)$ by 
Theorem \ref{MainTheoSolvGps}, we may also assume without loss that $P\o=Q\o$. 
But then $P=Q$ by Corollary \ref{CorQ1oQ2oGenCase}, a contradiction since $g\notin QR\o(G)$. 

$(b)$. 
Let $A$ be the pull back in $G$ of the large set of $G/R\o(G)$ provided in Theorem \ref{transfer2}(5), and 
$$A=A_1 \sqcup \cdots \sqcup A_{n(G)}$$ 
the pull back in $G$ of the corresponding partition of that large set equally provided in 
Theorem \ref{transfer2}(5). Here $n(G)$ is the number of ``generic types" of $G/R\o(G)$ 
computed with precision in Remark \ref{RemarknG}. By additivity of the dimension, 
$A$ is large in $G$ and each $A_i$ is weakly generic. 
Our claim is that for $Q$ a Cartan subgroup of $G$ and $a\in {Q\cap A_i}$ for some $i$, the set 
$[aQ\o]^{R\o(G)}$ is large in $aQ\o R\o(G)$. Since $Q\o$ normalizes the coset $aQ\o$, this is equivalent 
to showing that $[aQ\o]^{Q\o R\o(G)}$ is large in $aQ\o R\o(G)$. But by 
Theorem \ref{transfer2}(4-5) applied in $G/R\o(G)$, one can see that 
the largeness of the set of Cartan subgroups of $G$ and the additivity of the dimension 
forces $[aQ\o]^{Q\o R\o(G)}$ to be large in $aQ\o R\o(G)$. 
\qed

For instance, if $G$ is a definably connected real Lie group definable in an o-minimal expansion of $\R$, 
then its Cartan subgroups form a large subset by Fact \ref{Lie5}$(a)$ and the fact that 
density implies largeness for definable sets (as seen in the proof of Theorem \ref{transfer}). 
Hence, by Proposition \ref{PropLargeImpliesAll}, such a $G$ can produce a 
counterexample to neither Question \ref{QuestLiftCartans} nor Question \ref{QuestaQolqrgeGen}. 
Attacking Question \ref{QuestQlarges} in general would seem to rely on an abstract version of 
Fact \ref{Lie5}$(a)$, but with a priori no known abstract analog of regular elements (as in the proof 
of Fact \ref{Lie4}$(c)$) it seems difficult to find any spark plug. 

\section{Final remarks}\label{SectionFinalComments}

We begin this final section with additional comments on 
Question \ref{QuestLiftCartans} in special cases. 
If $G$ is a definably connected group definable in an o-minimal structure, 
then by Fact \ref{TheoCherlinZilberConjoMin} we have 
$$G/R(G)={G_1/R(G)} \times \cdots \times {G_n/R(G)}$$ 
for some definable subgroups $G_i$ containing $R(G)$ and such that $G_i/R(G)$ is definably simple. 
For each $i$, $G_i/R(G)$ is definably connected, and thus $G_i=G_i\o R(G)$. From the decomposition 
$G=G_1\cdots G_n$ we get $G=G_1\o \cdots G_n\o R(G)$. 
By definable connectedness of $G$ we also get a decomposition 
$$G=G_1\o \cdots G_n\o \leqno{(*)}$$
where each $G_i\o$ is definably connected, contains $R\o(G)$, and 
$G_i\o / R\o(G)$ 
is finite-by-(definably simple), as $R(G_i\o)={G_i\o \cap R(G)}$ and $G_i\o / R(G_i\o)$ is definably isomorphic 
to $G_i/R(G)$. We may analyze certain factors $G_i\o$ individually with the following. 

\begin{Fact}\label{FactGFiniteBySimpleAlg}
Let $\mathcal M$ be an o-minimal structure and $G$ a definably connected group definable in $\mathcal M$ with 
$R(G)=Z(G)$ finite and $G/R(G)$ definably simple. 
\begin{itemize}
\item[$(a)$] 
If $G/R(G)$ is stable as in the first case of Fact \ref{TheoCherlinZilberConjoMin}, 
then $G$ is (definably isomorphic in $\mathcal M$ to) an algebraic group over an algebraically 
closed field. 
\item[$(b)$] 
If $G/R(G)$ is definably compact, then $G$ is definably compact as well. 
\end{itemize}
\end{Fact}
\proof
As $G$ is definably connected and semisimple, there is an $\mathcal M$-definable real closed field $R$ such 
that $G$ is definably isomorphic in $\mathcal M$ to a semialgebraic group over the field of real algebraic 
numbers $R_{alg}\subseteq R$, by \cite[4.4(ii)]{MR2746030} or \cite{MR2765633}. 
In case $(a)$ our claim follows from \cite[6.3]{MR2746030} and thus we only have to consider case $(b)$. 
Assume towards a contradiction that $\alpha:(0,1)\rightarrow G$ is a continuous 
definable curve not converging in $G$. Since $G/Z(G)$ is definably compact, 
the composition of $\alpha$ with the projection $p:G\rightarrow G/Z(G)$ converges to a point 
$x\in G/Z(G)$. By \cite[Prop.2.11]{MR2114966}, $p$ is a definable covering map. 
In particular, there exists a definable open neighbourhood $U$ of $x$ in 
$G/Z(G)$ such that each definable connected component of $p^{-1}(U)$ is definably 
homeomorphic to $U$ via $p$. Since $\alpha$ does not converge to any point of 
$p^{-1}(x)=\{y_1,\cdots,y_{s}\}$, by o-minimality there exist definable open neighbourhoods 
$V_i\subseteq p^{-1}(U)$ of $y_i$ and $\delta\in (0,1)$ such that 
$\alpha(t)\notin {V_1\cup \cdots \cup V_{s}}$ for $t\in (\delta,1)$. Hence 
$p\circ \alpha(t)$ does not lie in the open neighbourhood $p(V_1)\cap \cdots \cap p(V_{s})$ 
of $x$ for $t\in (\delta,1)$, which is a contradiction. 
\qed

\begin{Corollary}\label{CorGGOOD}
If $G$ is as in Fact \ref{FactGFiniteBySimpleAlg}, case $(a)$ or $(b)$, then it has a single conjugacy class of 
Cartan subgroups, which are divisible and definably connected. 
\end{Corollary}
\proof
It is well known that in a connected reductive algebraic group over an algebraically closed field, 
Cartan subgroups are the selfcentralizing maximal algebraic tori, and are conjugate. 
They are isomorphic to a direct product of finitely many copies of the multiplicative group of the ground 
field (where the number of copies is the Lie rank of the group seen as a pure algebraic group). 
In particular they are divisible, and thus with no proper subgroup of finite index. 
In the definably compact case we refer to Corollary \ref{CorStructCompactGps}, 
getting the divisibility from the definable connectedness in this case. 
\qed

Consider the decomposition $(*)$ of a definably connected group $G$ as above, and let $I={\{1,\cdots , n\}}$. 
Let $I_1$ be the subset of elements $i\in I$ such that $G_i\o/R\o(G_i\o)$ is stable (as a pure group) or 
definably compact. Notice that, by Fact \ref{FactGFiniteBySimpleAlg}, it suffices to require the definably simple group 
$G_i\o/R(G_i\o)$ to be stable (as a pure group) or definably compact. 
Let $I_2$ be the subset of elements $i\in I$ such that Cartan subgroups of $G_i\o/R\o(G_i\o)$ 
are definably connected. Finally, let $I_3$ be the subset of elements $i\in I$ such that in $G_i\o$ 
Question \ref{QuestLiftCartans} has a positive answer for any Cartan subgroup. 
Corollary \ref{CorGGOOD} shows that 
$I_1 \subseteq I_2$ and Lemma \ref{LiftCarterGeneral} shows that $I_2\subseteq I_3$. Hence 
$$I_1\subseteq I_2 \subseteq I_3 \subseteq I$$
and the inclusion $I_1\subseteq I_3$ reads informally as the fact that the definably simple factors of $G/R(G)$ 
which are algebraic or compact cannot produce any counterexample to the lifting problem 
of Question \ref{QuestLiftCartans}. More precisely, we have the following statement. 

\begin{Remark}
If $I_2=I$, then $G$ cannot produce any counterexample to the lifting problem 
of Question \ref{QuestLiftCartans}. 
\end{Remark}
\proof
First one can check that, modulo $R\o(G)$, the decomposition $(*)$ of $G$ becomes a central product: 
$$G/R\o(G)={G_1\o/R\o(G) \ast \cdots \ast G_1\o/R\o(G)}.$$ 
Indeed, if $i\neq j$, then 
$[G_i\o, G_j\o] \leq R(G)$, and $R(G)$ is finite modulo $R\o(G)$. Hence any element in 
$G_i\o/R\o(G)$ has a centralizer of finite index in the other factor $G_j\o/R\o(G)$, which must 
then be the full factor $G_j\o/R\o(G)$ by definable connectedness. Therefore the factors $G_i\o/R\o(G)$ pairwise 
commute, as claimed. 
Now Lemma \ref{LemElemCartanCentProd} gives that Cartan subgroups of $G/R\o(G)$ are exactly of the 
form ${Q_1/R\o(G) \ast \cdots \ast Q_n/R\o(G)}$ with, for each $i$, $Q_i/R\o(G)$ a Cartan subgroup of 
${G_i\o/R\o(G)}$. 

Assuming now that $I_2=I$ we get that, for each $i$, each Cartan subgroup $Q_i/R\o(G)$ of $G_i\o/R\o(G)$ 
is definably connected. We then see that Cartan subgroups of $G/R\o(G)$ must be definably connected 
as well. Now Lemma \ref{LiftCarterGeneral} implies that Question \ref{QuestLiftCartans} is positively 
satisfied for every Cartan subgroup of $G$ 
(and that such Cartan subgroups of $G$ are all definably connected and Carter subgroups 
by Corollary \ref{GConResCartan=Carter}). 
\qed

The decomposition $(*)$ of a definably connected group $G$ as above is also convenient for 
describing the various relative Weyl groups. If $Q$ is a Cartan subgroup of $G$, then we still have that 
$N_G(Q\o)=N_G(Q)$ by Corollary \ref{CorQ1oQ2oGenCase}. If Question \ref{QuestLiftCartans} is positively 
satisfied for $Q$, then retaining the 
notation of Section \ref{SectionGeneralCase} and using the notation ``$\overline{\phantom{H}}$" for 
quotients modulo $R\o(G)$ we get, as after Lemma \ref{LemKQEqualsPRQR}, that 
$$W(\overline{G},\overline{K_Q}) \simeq N_G(Q)/Q.$$
We also see, with Theorem \ref{TheoGen(2)} or just Lemma \ref{LemElemCartanModZn}$(a)$, that 
$\overline{R(G)}$ does not contribute to the Weyl group $W(\overline{G},\overline{K_Q})$. 
Hence the latter is isomorphic to the direct product of the Weyl groups in $G_i/R(G)$ 
relative to the factors of $QR(G)/R(G)$ in its 
decomposition along the decomposition ${G_1/R(G)} \times \cdots \times {G_n/R(G)}$ of 
$G/R(G)$ (Corollary \ref{LemElemCartanDirProd}). 
Since the group $N_G(Q)/Q$ is isomorphic to $W(\overline{G},\overline{K_Q})$, it has the 
same isomorphism type and may be called the {\em Weyl group relative to} $Q$. 

Without assuming the exact lifting of Question \ref{QuestLiftCartans} for the Cartan subgroup $Q$ 
we only get, with Corollary \ref{CorTheoGen(2)} and 
as after Lemma \ref{LemKQEqualsPRQR}, that 
$$W(\overline{G},\overline{K_Q}) \simeq (N_G(Q)/Q)/({\hat{Q}/Q}).$$ 
In this case the Weyl group $W(\overline{G},\overline{K_Q})$ has the same description as above, but 
$N_G(Q)/Q$ just has a quotient isomorphic to $W(\overline{G},\overline{K_Q})$. 

We finish on a more model-theoretic note. 

\begin{Proposition}\label{PropFinal}
Let $\mathcal{M}$ be an o-minimal structure, $A\subseteq M$ 
a set of parameters such that ${\hbox{\rm dcl}}_{\mathcal M}(A)\preceq {\mathcal M}$, and $G$ a group definable in 
$\mathcal{M}$ over $A$. Then $G$ has a finite system of representatives of Cartan 
(resp. Carter) subgroups, both definable over $A$. 
\end{Proposition}

Of course, having now Theorem \ref{MainTheo} at hand, the proof of Proposition \ref{PropFinal} 
is the same as in 
Lemma \ref{lemparam} and Corollary \ref{corparam}. As seen in the proof of 
Lemma \ref{lemparam}, when $\mathcal M$ expands an ordered abelian group, examples of $A$ 
such that ${\hbox{\rm dcl}}_{\mathcal M}(A)\preceq {\mathcal M}$ include any $A$ not contained in $\{0\}$. 

\bibliographystyle{plain}
\bibliography{biblio}

\end{document}